# A Bernstein-Von Mises Theorem for discrete probability distributions

## S. Boucheron[*]

LPMA, CNRS *and Université Paris-Diderot*
*e-mail:* stephane.boucheron@math.jussieu.fr
*url:* http://www.proba.jussieu.fr/~boucheron

## E. Gassiat[†]

CNRS *and Université Paris-Sud 11,*
*e-mail:* elisabeth.gassiat@math.u-psud.fr
*url:* http://www.math.u-psud.fr/~gassiat

**Abstract:** We investigate the asymptotic normality of the posterior distribution in the discrete setting, when model dimension increases with sample size. We consider a probability mass function $\theta_0$ on $\mathbb{N} \setminus \{0\}$ and a sequence of truncation levels $(k_n)_n$ satisfying $k_n^3 \leq n \inf_{i \leq k_n} \theta_0(i)$. Let $\hat{\theta}$ denote the maximum likelihood estimate of $(\theta_0(i))_{i \leq k_n}$ and let $\Delta_n(\theta_0)$ denote the $k_n$-dimensional vector which $i$-th coordinate is defined by $\sqrt{n}(\hat{\theta}_n(i) - \theta_0(i))$ for $1 \leq i \leq k_n$. We check that under mild conditions on $\theta_0$ and on the sequence of prior probabilities on the $k_n$-dimensional simplices, after centering and rescaling, the variation distance between the posterior distribution recentered around $\hat{\theta}_n$ and rescaled by $\sqrt{n}$ and the $k_n$-dimensional Gaussian distribution $\mathcal{N}(\Delta_n(\theta_0), I^{-1}(\theta_0))$ converges in probability to 0. This theorem can be used to prove the asymptotic normality of Bayesian estimators of Shannon and Rényi entropies.

The proofs are based on concentration inequalities for centered and non-centered Chi-square (Pearson) statistics. The latter allow to establish posterior concentration rates with respect to Fisher distance rather than with respect to the Hellinger distance as it is commonplace in non-parametric Bayesian statistics.

Received July 2008.

## Contents

---

[*]Supported by ANR Project TAMIS.
[†]Supported by NOE PASCAL2

## 1. Introduction

The classical Bernstein-Von Mises Theorem asserts that for regular (Hellinger differentiable) parametric models, under mild smoothness conditions on the prior distribution, after centering around the maximum likelihood estimate and rescaling, the posterior distribution of the parameter is asymptotically Gaussian and that the limiting covariance matrix coincides with the inverse of the Fisher information matrix. This theorem provides a frequentist perspective on the Bayesian methodology and elements for reconciliation of the two approaches. In regular parametric models, Bernstein-von Mises theorems motivate the interchange of Bayesian credible sets and frequentist confidence regions. Refinements of the Bernstein-von Mises theorem have also proved helpful when analyzing the redundancy of universal coding for smoothly parametrized classes of sources over finite alphabets.

The proof of the classical Bernstein-Von Mises theorem relies on rather sophisticated arguments. Some of them seem to be tied up with the finite dimensionality of the considered models. Hence, extensions of Bernstein-von Mises theorems to non-parametric and semi-parametric settings have both received deserved attention and shown moderate progress during the last four decades. Soon after Bayesian inference was put on firm frequentist foundations by Doob (1949), Schwartz (1965) and others, Freedman (1963) (see also Freedman, 1965) pointed out that even when dealing with the simplest possible case, that of independent, identically distributed, discrete observations, there is no such thing as a general posterior consistency result let alone a general Bernstein-Von Mises Theorem. Moreover, according to the evidence presented by Freedman (1965), it is mandatory to focus moderately large classes of distributions. Despite such early negative results, non-parametric Bayesian theory has been progressing at a steady pace. The framework of empirical process theory has enabled to provide sufficient conditions for posterior consistency and to relate posterior concentration rates to model complexity (Ghosal and van der Vaart, 2007b, 2001; Ghosal et al., 2000).

Among the different approaches to non-parametric inference, using simple models with increasing dimensions has attracted attention in the context of maximum likelihood inference (Portnoy, 1988; Fan and Truong, 1993; Fan et al., 2001; Fan, 1993) and in the context of Bayesian inference (Ghosal, 2000). The last reference is especially relevant to this paper. Therein, S. Ghosal considers nested sequences of exponential models satisfying a number of assumptions involving the growth rate of models with sample size, the growth rate of the determinant of the Fisher information matrix with respect to model dimension (and thus sample size), prior smoothness, and moment bounds for score functions in small Kullback-Leibler balls located around the sampling probability (those conditions will be explained and compared with our own conditions in Section 3.1). S. Ghosal proves a Bernstein-Von Mises Theorem (Ghosal, 2000, Theorem 2.3) for the log-odds parametrization, partially building on previous results from Portnoy (1988) concerning maximum likelihood estimates. However our objectives significantly differ from those of S. Ghosal. In (Ghosal, 2000), the main application of non-parametric Bernstein-Von Mises Theorems for multinomial models seems to be non-parametric density estimation using histograms. This framework justifies special attention to multinomial distributions which are almost uniform. Our ultimate goal is quite different. In information-theoretical language, we are interested in investigating memoryless sources over infinite alphabets as in (See Kieffer, 1978; Gyorfi et al., 1993; Boucheron et al., 2009, and references therein). In Information Theory, refinements of Bernstein-Von Mises Theorems allow to investigate the so-called maximin redundancy of universal coding over parametric classes of sources (Clarke and Barron, 1994). In Information Theory, a source over a (countable alphabet) is a probability distribution over the set of infinite sequences of symbols from the alphabet. The redundancy of a (coding) probability distribution with respect to a source on a given (finite) sequence of symbols is the logarithm of the ratio between the probability of the sequence under the source and under the coding probability. In universal coding theory, average redundancy with respect to a prior distribution over sources can be written as the difference between the (differential) Shannon entropy of the prior distribution and the average value of the (differential) entropy of the conditional posterior distribution. Thanks to non-trivial refinements of the Bernstein-Von Mises Theorem, the latter conditional entropy can be approximated by the (differential) entropy of a Gaussian distribution which covariance matrix is the inverse of the Fisher information matrix defined by the source under consideration. This elegant approach provides sharp asymptotic and non-asymptotic results when dealing with classes of sources which are soundly parameterized by subsets of finite-dimensional spaces (See Clarke and Barron, 1990, for precise definitions). When turning to larger classes of sources, for example toward memoryless sources over countable alphabets (Boucheron et al., 2009), this approach to the characterization of maximin redundancy has not (yet) been carried out. A major impediment is the current unavailability of adequate non-parametric Bernstein-Von Mises Theorems.

This paper is a first step in developing the Bayesian tools that are useful to precisely quantify the minimax redundancy of universal coding of non-parametric classes of sources over infinite alphabets. Because of our ultimate goals, we cannot focus on almost uniform multinomial models. We are specifically interested in situations where the sampling probability mass functions decay at a prescribed rate (say algebraic or exponential) as in (Boucheron et al., 2009).

As pointed out by Ghosal, in models with an increasing number of parameters, justifying the asymptotic normality of the posterior distribution is more involved, and precisely characterizing under which conditions on prior and sampling distribution this asymptotic normality holds remains an open-ended question. For example, in the context of discrete distributions, several ways of defining the divergence between distributions look reasonable. Most of the recent work on non-parametric Bayesian statistics dealt with posterior concentration rates and has been developed using Hellinger distance (Ghosal et al., 2000; Ghosal and van der Vaart, 2007b, 2001). One may wonder whether some posterior concentration rate results obtained using Hellinger metrization can be strengthened. It is not clear how to tackle this issue in full generality. In this paper, taking advantage of the peculiarities of our models, we use another, demonstrably stronger, information divergence, the Fisher ($\chi^2$) "distance" and establish posterior concentration rates with respect to Fisher balls (see 3.6). The proof relies on known concentration inequalities for centered $\chi^2$ (Pearson) statistics and (apparently) new concentration inequalities for non-centered $\chi^2$ statistics.

Paraphrasing van der Vaart (1998), as the notion of convergence in the Bernstein-Von Mises Theorem is a rather complicated one, the expected reward, once such a Theorem has been proved, is that "nice" functionals applied to the posterior laws should converge in distribution in the usual sense. An obvious candidate for deriving that kind of method is a Bayesian variation on the Delta method. However, we are facing here two kinds of obstacles. On the one hand, we cannot rely on the availability of a Bernstein-Von Mises Theorem when considering the infinite-dimensional model (Freedman, 1963, 1965). This precludes using the traditional functional Delta method as described for example in (van der Vaart and Wellner, 1996; van der Vaart, 1998). On the other hand, when considering models of increasing dimensions, a variant of the Delta method has to be derived in an *ad hoc* manner. This is what we do. We assess this rule of thumb by examining plug-in estimates of Shannon and Rényi entropies. Such functionals characterize the compressibility of a given probability distribution (Csiszár and Körner, 1981; Cover and Thomas, 1991; Gallager, 1968). The problem of estimating such functionals has been investigated by Antos and Kontoyiannis (2001) and Paninski (2004). It has been checked there that plug-in estimates of the Shannon and Rényi entropies are consistent and some lower and upper bounds on the rate of convergence have been proposed. Up to our knowledge, classes of distributions for which plug-in estimates satisfy a central limit theorem have not been systematically characterized. Here, the Bernstein-Von Mises Theorem allows to derive central limit theorems for Bayesian entropy estimators (see Theorem 3.12) and provides the basis for constructing Bayesian credible sets. In the present context, those credible sets are known to coincide asymptotically with Bayesian bootstrap confidence regions (Rubin, 1981).

The paper is organized as follows. In Section 2, the framework and notation of the paper are introduced. A few technical conditions warranting local asymptotic normality when handling models of increasing dimensions are also stated. The main results of the paper are presented in Section 3. The nonparametric Bernstein-Von Mises Theorem (3.7) is described in Subsection 3.1. It is complemented by a posterior concentration lemma (3.6) that might be interesting in its own right. A roadmap of the proof of the Bernstein-Von Mises Theorem is stated thereafter. In Paragraph 3.2, the asymptotic normality of Bayesian estimators of various entropies is derived using the non-parametric Bernstein-Von Mises Theorem and various tail bounds for quadratic forms that are also useful in the derivation of the Bernstein-Von Mises theorem. In Paragraph 3.3, sequences of Dirichlet priors are checked to satisfy the conditions of the Bernstein-Von Mises Theorem. The main results of the paper are illustrated on the envelope classes investigated by Boucheron et al. (2009). In Subsection 3.5, the setting of Theorem 3.7 is compared with the framework described in (Ghosal, 2000). In Subsection 3.6, the posterior concentration lemma is compared with related recent results in non-parametric Bayesian statistics. The Proof of the Bernstein-Von Mises Theorem is given in Section 4. It adapts Le Cam's proof (Le Cam and Yang, 2000; van der Vaart, 2002) to the non-parametric setting using a collection of old and new non-asymptotic tail bounds for chi-square statistics. The proof of the asymptotic normality of Bayesian entropy estimators is given in Section 5. It relies on the Bernstein-Von Mises Theorem and on the aforementioned tail bounds for chi-square statistics.

## 2. Notation and background

This section describes the statistical framework we will work with, as well as the behavior of likelihood ratios in this framework. At the end of the section, a useful contiguity result is stated.

Throughout the paper, $\theta = (\theta(i))_{i \in \mathbb{N}_*}$ denotes a probability mass function over $\mathbb{N}_* = \mathbb{N} \setminus \{0\}$ and $\Theta$ denotes the set of probability mass functions over $\mathbb{N}_*$. If the sequence $\mathbf{x} = x_1, \ldots, x_n$ denotes a sample of $n$ elements from $\mathbb{N}_*$, let $N_i$ denote the number of occurrences of $i$ in $\mathbf{x}$: $N_i(\mathbf{x}) = \sum_{j=1}^{n} 1_{\mathbf{x}_j = i}$. The log-likelihood function maps $\Theta \times \mathbb{N}_*^n$ toward $\mathbb{R}$:

$$\ell_n(\theta, \mathbf{x}) = \sum_{i \geq 1} N_i \log \theta(i).$$

When the sample $\mathbf{x}$ is clear from context, $\ell_n(\theta, \mathbf{x})$ is abbreviated into $\ell_n(\theta)$.

Throughout the paper, $\theta_0$ denotes the (unknown) probability mass function under which samples are collected. Let $\Omega = \mathbb{N}_*^{\mathbb{N}}$, let $X_1, \ldots, X_n, \ldots$ denote the coordinate projections. Then $\mathbb{P}_0$ denotes the probability distribution over $\Omega$ (equipped with the cylinder $\sigma$-algebra $\mathcal{F}$), satisfying

$$\mathbb{P}_0\{\wedge_{i=1}^n X_i = x_i\} = \prod_{i=1}^n \theta_0(x_i).$$

Recall that the maximum likelihood estimator $\widehat{\theta}$ of $\theta_0$ on a sample $\mathbf{x}$ is given by the empirical probability mass function: $\widehat{\theta}(i) = N_i/n$.

Let $k$ denote a positive integer that may and should depend on the sample size $n$. We will be interested in the estimation of the $\theta_0(i)$ for $i = 1, \ldots, k$. In this respect, all the useful information is conveyed by the counts $N_i$, $i = 1, \ldots, k$, or equivalently in what will be called the truncated version of the sample. The truncated version of sample $\mathbf{x}$ is denoted by $\tilde{\mathbf{x}}$ and constructed as follows

$$\tilde{\mathbf{x}}_i = \begin{cases} \mathbf{x}_i & \text{if } \mathbf{x}_i \leq k \\ 0 & \text{otherwise}. \end{cases}$$

The counter $N_0$ is defined as the number of occurrences of 0 in $\tilde{\mathbf{x}}$: $N_0(\mathbf{x}) = \sum_{i>k} N_i(\mathbf{x})$. The image of $\theta \in \Theta$ by truncation is a p.m.f. over $\{0, \ldots, k\}$, it is still denoted by $\theta$ with $\theta(0) = \sum_{i>k} \theta(i)$. Let $\Theta_k$ denote the set of p.m.f. over $\{0, \ldots, k\}$. In the sequel, depending on context, $\theta_0$ may denote either the p.m.f. on $\mathbb{N}_*$ from which the sample is drawn or its image by truncation at level $k$.

Henceforth, $\theta \in \Theta_{k_n}$ may denote either $(\theta(i))_{0 \leq i \leq k_n}$ or its projection on the $k_n$ last coordinates $(\theta(i))_{1 \leq i \leq k_n}$; in the same way, if $h$ denotes a vector $(h(i))_{0 \leq i \leq k_n}$ in $\mathbb{R}^{k_n+1}$ such that $\sum_{i=0}^{k_n} h(i) = 0$, $h$ may also denote its projection on the $k_n$ last coordinates $(h(i))_{1 \leq i \leq k_n}$ depending on the context.

For a given sample $\mathbf{x}$, the score function is the gradient of the log-likelihood at $\theta \in \Theta_k$, for $i \in \{1, \ldots, k\}$: $(\dot{\ell}_n(\theta))_i = N_i/\theta(i) - N_0/\theta(0)$. Assume all components of $\theta \in \Theta_k$ are positive, then the information matrix $I(\theta)$ is defined as

$$I(\theta) = \frac{1}{n} \mathbb{E}_\theta \left[ \dot{\ell}_n(\theta) \dot{\ell}_n^T(\theta) \right] = \text{Diag}\left(\frac{1}{\theta(i)}\right)_{1 \leq i \leq k} + \frac{1}{\theta(0)} \mathbf{1}\mathbf{1}^T$$

and its inverse is

$$I^{-1}(\theta) = \text{Diag}(\theta(i))_{1 \leq i \leq k} - \begin{pmatrix} \theta(1) \\ \vdots \\ \theta(k) \end{pmatrix} \begin{pmatrix} \theta(1) & \ldots & \theta(k) \end{pmatrix}.$$

It can be checked that $\det(I(\theta)) = \prod_{i=0}^k \theta^{-1}(i)$. The pseudo-sufficient statistic $\Delta_n(\theta)$ is defined as

$$\Delta_n(\theta) = \frac{1}{\sqrt{n}} I^{-1}(\theta) \dot{\ell}_n(\theta) = \sqrt{n}(\hat{\theta} - \theta).$$

Note that $\sqrt{n}(\widehat{\theta} - \theta_0) = \Delta_n(\theta_0)$ and that this $k$-dimensional random vector has covariance matrix $I^{-1}(\theta_0)$. Moreover for each positive $\theta \in \Theta_k$, $\Delta_n^T(\theta)I(\theta)\Delta_n(\theta)$ coincides with the Pearson $\chi^2$ statistics.

$$\Delta_n^T(\theta)I(\theta)\Delta_n(\theta) = \sum_{i=0}^{k} \frac{(N_i - n\theta(i))^2}{n\theta(i)}.$$

Let $k_n$ denote a truncation level. If $h$ belongs to $\mathbb{R}^{k_n+1}$ and satisfies $\sum_{i=0}^{k_n} h(i) = 0$, let $\sigma_n(h)$ be defined by

$$\sigma_n^2(h) = \sum_{i=0}^{k_n} \frac{h^2(i)}{\theta_0(i)} = h^T I(\theta_0) h,$$

where we agree on the following convention: if $\theta_0(i) = 0$ and $h(i) = 0$, then $h^2(i)/\theta_0(i) = 0$. The set $\mathcal{E}_{\theta_0,k_n}(M)$ is the intersection of a $k_n$-dimensional subspace with an ellipsoid in $\mathbb{R}^{k_n+1}$.

$$\mathcal{E}_{\theta_0,k_n}(M) = \left\{ h \,:\, \sigma_n^2(h) \leq M, \sum_{i=0}^{k_n} h(i) = 0, \, h(i) \geq -\sqrt{n}\theta_0(i), i = 0, \ldots, k_n \right\}.$$

In the parametric setting, that is when $k_n$ remains fixed, Le Cam's proof of the Bernstein-Von Mises Theorem (van der Vaart, 1998; van der Vaart, 2002) is made significantly more transparent by resorting to a contiguity argument. In order to adapt this argument to our setting, we need to formulate two conditions.

In the sequel $(k_n)_{n \in \mathbb{N}}$ denotes a non-decreasing sequence of truncation levels.

**Condition 2.1.** *The p.m.f. $\theta_0$ and the sequence $(k_n)_{n \in \mathbb{N}}$ satisfy*

$$n \inf_{i \leq k_n} \theta_0(i) \to +\infty.$$

Let $(h_n)_{n \in \mathbb{N}}$ denote a sequence of elements from $\mathbb{R}^{k_n+1}$ such that for each $n$, $\sum_{i=0}^{k_n} h_n(i) = 0$. The sequence $(h_n)_{n \in \mathbb{N}}$ is said to be tangent at the p.m.f. $\theta_0$ if the following condition is satisfied.

**Condition 2.2.** *There exists a positive real $\sigma$ such that the sequence $\sigma_n^2(h_n)$ tends toward $\sigma^2 > 0$.*

The probability distribution $\mathbb{P}_{n,h}$ over $\{0, \ldots, k_n\}^n$ is the product distribution defined by the perturbed p.m.f. $\theta_0(i) + h(i)/\sqrt{n}$ if $0 < \theta_0(i) + \frac{h(i)}{\sqrt{n}} < 1$ for all $i$ in $\{0, \ldots, k_n\}$. We are now equipped to state the building block of the contiguity argument: the proof is given in the appendix (A).

**Lemma 2.3.** *Let $\theta_0$ denote a probability mass function over $\mathbb{N}_*$. If the sequence of truncation levels $(k_n)_{n \in \mathbb{N}}$ satisfy Condition 2.1 and if the sequence $(h_n)_{n \in \mathbb{N}}$ satisfies the tangency Condition 2.2 then the sequences $(\mathbb{P}_{n,h})_n$ and $(\mathbb{P}_{n,0})_n$ are mutually contiguous, that is, for any sequence $(B_n)$ of events where for each $n$, $B_n \subseteq \{0, \ldots, k_n\}^n$, the following holds:*

$$\lim_n \mathbb{P}_{n,h}\{B_n\} = 0 \Leftrightarrow \lim_n \mathbb{P}_{n,0}\{B_n\} = 0.$$

Note that throughout the paper, we use De Finetti's convention: if $(\Omega, \mathcal{F}, P)$ denotes a probability space, $Z$ a random variables on $\Omega, \mathcal{F}$, then $PZ = P[Z] = P(Z)$ denotes the expected value of $Z$ (provided it is well-defined, that is $PZ_+$ and $PZ_-$ are not both infinite. If $A$ denotes an event, then $P\{A\} = P\mathbf{1}_A$.

## 3. Main results

In a Bayesian setting, the set of parameters is endowed with a prior distribution. In this paper, we consider a sequence of prior distributions $(W_n)_{n\in\mathbb{N}}$ matching the non-decreasing sequence of truncation levels we use. Let $W_n$ be a prior probability distribution for $(\theta(i))_{1\leq i\leq k_n}$ such that $\theta = (\theta(i))_{0\leq i\leq k_n} \in \Theta_{k_n}$. Henceforth, we assume that $W_n$ has a density $w_n$ with respect to Lebesgue measure on $\mathbb{R}^{k_n}$. Let $T = (\tau(i))_{0\leq i\leq k_n}$ be a random variable such that $(\tau(i))_{1\leq i\leq k_n}$ is distributed according to $W_n$ and $\tau(0) = 1 - \sum_{i=1}^{k_n} \tau(i)$. Conditionally on $T = \theta$, $(X_n)_{n\in\mathbb{N}}$ is a sequence of independent random variables distributed according to the p.m.f. $\theta$.

### 3.1. Non parametric Bernstein-Von Mises Theorem

Let $H_n$ be the random variable $H_n = \sqrt{n}(\tau(i) - \theta_0(i))_{1\leq i\leq k_n}$, and $P_{H_n|X_{1:n}}$ its posterior distribution, that is its distribution conditionally to the observations $X_{1:n} = (X_1, \ldots, X_n)$. If the truncation level $k_n = k$ (that is the dimension of the parameter space $\Theta_{k_n}$) is a constant integer, the classical parametric Bernstein-Von Mises Theorem asserts that the sequence of posterior distributions is asymptotically Gaussian with centerings $\Delta_n(\theta_0) = \sqrt{n}(\widehat{\theta} - \theta_0)$ and variance $I^{-1}(\theta_0)$ if the observations $X_{1:n}$ are independently distributed according to $\theta_0$.

Theorem 3.7 below asserts that under adequate conditions on the sequence of priors $W_n$ and on the tail behavior of $\theta_0$, the Bernstein-Von Mises Theorem still holds provided the truncation levels $k_n$ do not increase too fast toward infinity.

For any sequence of prior distributions $(W_n)_n$, for a sequence $M_n$ of real numbers increasing to $+\infty$, and a sequence $(k_n)_n$ of truncation levels that satisfy Condition 2.1, we will use the following three conditions in order to establish the three propositions the Bernstein-Von Mises Theorem depends on.

**Condition 3.1.** *The sequence of truncation levels $(k_n)_n$ and radii $M_n$ satisfies*

$$M_n = o\left(\left(n \inf_{i\leq k_n} \theta_0(i)\right)^{1/3}\right), \tag{3.2}$$

$$k_n = o(M_n). \tag{3.3}$$

Requiring a prior smoothness condition is commonplace when establishing asymptotic normality of posterior distribution in parametric settings.

**Condition 3.4.** (PRIOR SMOOTHNESS)

$$\sup_{h,g\in\mathcal{E}_{\theta_0,k_n}(M_n)} \frac{w_n\left(\theta_0 + \frac{h}{\sqrt{n}}\right)}{w_n\left(\theta_0 + \frac{g}{\sqrt{n}}\right)} \to 1,$$

Requiring a prior concentration condition, sometimes called a small ball probability conditions is usual in non-parametric Bayesian statistics.

**Condition 3.5.** (PRIOR CONCENTRATION)

$$\frac{k_n}{2} \log(n) \vee \log(\det(I(\theta_0))) \vee -\log(w_n(\theta_0)) = o(M_n)$$

where $\det(I(\theta_0)) = \prod_{i=0}^{k_n} \theta_0^{-1}(i)$.

Note that the prior concentration condition entails the second condition in Condition 3.1.

The next lemma which is proved in Section 4.3 asserts that under mild conditions, the posterior distribution concentrates on $\chi^2$ (Fisher) balls centered around maximum likelihood estimates.

**Lemma 3.6.** (POSTERIOR CONCENTRATION) *If the p.m.f. $\theta_0$ and the sequence of truncation levels $(k_n)_n$ both satisfy Conditions (2.1, 3.4, 3.5) and if $M_n = o(n \inf_{i \leq k_n} \theta_0(i))$ then under $\mathbb{P}_{n,0}$*

$$\mathbb{P}_{n,0} P_{H_n|X_{1:n}} \left( H_n^T I(\theta_0) H_n \geq M_n \right) = \mathbb{P}_{n,0} P_{H_n|X_{1:n}} \left( H_n \notin \mathcal{E}_{0,k_n}(M_n) \right) \to 0.$$

This posterior concentration lemma allows to recover the parametric posterior concentration phenomenon if truncation levels remain fixed and strengthens the generic non-parametric posterior concentration theorem from Ghosal et al. (2000).

**Theorem 3.7.** (A NON-PARAMETRIC BERNSTEIN-VON MISES THEOREM) *If the sequence of truncation levels $(k_n)_{n \in \mathbb{N}}$, $k_n \to +\infty$, and the p.m.f. over $\mathbb{N}_*$, $\theta_0$ satisfy Condition 2.1, and if there is an increasing sequence $(M_n)_n$ tending to infinity such that 3.1, 3.4 and 3.5 hold, then*

$$\mathbb{P}_{n,0} \left\| \mathcal{N}_{k_n}(\Delta_n(\theta_0), I^{-1}(\theta_0)) - P_{H_n|X_{1:n}} \right\| \to 0$$

*where $\| \cdot \|$ denotes the total variation norm.*

A comparison of the Theorem with respect to previous results available in the literature (Ghosal and van der Vaart, 2007b,a; Ghosal et al., 2000; Ghosal, 2000) is given at the end of the Section.

**Remark 3.8.** A corollary of the Bernstein-Von Mises Theorem is that

$$\mathbb{P}_{n,0} \left( P_{H_n|X_{1:n}} \left\{ H_n^T I(\theta_0) H_n \geq u_n \right\} \right) \to 0$$

if and only if $u_n/k_n \to \infty$.

The proof of Theorem 3.7 is organized along the lines of Le Cam's proof of the parametric Bernstein-Von Mises Theorem as exposed by A. van der Vaart in (van der Vaart, 1998) (see also van der Vaart (2002)).

*Roadmap of the proof of the Bernstein Von-Mises theorem.* If $P$ is any probability distribution on $\mathbb{R}^{k_n}$ and $M > 0$ is any positive real, let $P^M$ be the conditional probability distribution on the ellipsoid $\{u \in \mathbb{R}^{k_n} : u^T I(\theta_0) u = \sigma_n^2(u) \leq M\}$. For any measurable set $B$,

$$P^M(B) = \frac{P\{B \cap \{u \,:\, u^T I(\theta_0) u \leq M\}\}}{P\{u \,:\, u^T I(\theta_0) u \leq M\}}.$$

To alleviate notations, we will use the shorthands $\mathcal{N}_{k_n}$ and $\mathcal{N}_{k_n}^{M_n}$ to denote the (random) distributions $\mathcal{N}_{k_n}(\Delta_n(\theta_0), I^{-1}(\theta_0))$ and $\mathcal{N}_{k_n}^{M_n}(\Delta_n(\theta_0), I^{-1}(\theta_0))$.

From the triangle inequality, if follows that:

$$\begin{aligned}
&\left\|\mathcal{N}_{k_n}(\Delta_n(\theta_0), I^{-1}(\theta_0)) - P_{H_n|X_{1:n}}\right\| \\
&\leq \quad \left\|\mathcal{N}_{k_n} - \mathcal{N}_{k_n}^{M_n}\right\| + \left\|\mathcal{N}_{k_n}^{M_n} - P_{H_n|X_{1:n}}^{M_n}\right\| + \left\|P_{H_n|X_{1:n}}^{M_n} - P_{H_n|X_{1:n}}\right\|.
\end{aligned}$$

The proof of Theorem 3.7 boils down to checking that each of the three terms on the right-hand side tends to 0 in $\mathbb{P}_{n,0}$ probability.

The first term avers to be the easiest to control thanks to the well-known concentration properties of the Gaussian distribution. Upper bounding the middle term is arguably the most delicate part of the proof. The posterior concentration Lemma allows to deal with the third term.

Let us call $\textsc{nv}(M_n)$ the middle term

$$\left\|\mathcal{N}_{k_n}^{M_n} - P_{H_n|X_{1:n}}^{M_n}\right\|.$$

The posterior density is proportional to the product of the prior density and of the likelihood function. Hence, controlling the variation distance between $\mathcal{N}_{k_n}^{M_n}$ and $P_{H_n|X_{1:n}}^{M_n}$ requires a good understanding of log-likelihood ratios. A quadratic Taylor expansion of the log-likelihood ratio leads to:

$$\begin{aligned}
\log \frac{\mathbb{P}_{n,h}}{\mathbb{P}_{n,0}}(\mathbf{x}) &= \sum_{i=0}^{k_n} N_i \log\left(1 + \frac{h(i)}{\sqrt{n}\theta_0(i)}\right) \\
&= Z_n(h) - \frac{1}{2n}\sum_{i=0}^{k_n} N_i \frac{h^2(i)}{\theta_0^2(i)} + \frac{1}{n}\sum_{i=0}^{k_n} N_i \frac{h^2(i)}{\theta_0^2(i)} R\left(\frac{h(i)}{\sqrt{n}\theta_0(i)}\right) \\
&= Z_n(h) - \frac{\sigma_n^2(h)}{2} + \frac{A_n(h)}{2} + C_n(h)
\end{aligned}$$

where $R(u) = \frac{1}{u^2}(\log(1+u) - u - \frac{u^2}{2})$ satisfies $R(u) = O(u)$ as $u$ tends toward 0 and

$$Z_n(h) = \frac{1}{\sqrt{n}}\sum_{i=0}^{k_n} N_i \frac{h(i)}{\theta_0(i)},$$

$$A_n(h) = \sigma_n^2(h) - \frac{1}{n}\sum_{i=1}^{k_n} N_i \frac{h(i)^2}{\theta_0(i)^2}$$

and
$$C_n(h) = \frac{1}{n} \sum_{i=1}^{k_n} N_i \frac{h(i)^2}{\theta_0(i)^2} R\left(\frac{h(i)}{\sqrt{n}\theta_0(i)}\right).$$

Performing algebra along the lines described in (van der Vaart, 2002, P. 142) (computational details are given in the Appendix, see Section C), leads to

$$\mathrm{NV}(M_n) \leq$$
$$\iint \left(1 - \frac{w_n\left(\theta_0 + \frac{g}{\sqrt{n}}\right)}{w_n\left(\theta_0 + \frac{h}{\sqrt{n}}\right)} e^{\frac{A_n(g)-A_n(h)}{2}+C_n(g)-C_n(h)}\right)^+ \mathrm{d}\mathcal{N}_{k_n}^{M_n}(g) \mathrm{d}P_{H_n|X_{1:n}}^{M_n}(h). \tag{3.9}$$

We prove in Section 4.1 that the decay of $\mathrm{NV}(M_n)$ depends on prior smoothness around $\theta_0$ and on the ratio between $M_n$ and $(n \inf_{i \leq k_n} \theta_0(i))^{1/3}$:

**Proposition 3.10.**

$$\mathbb{P}_{n,0}\left(\mathrm{NV}(M_n)\right) = O\left(\sqrt{\frac{M_n^3}{n \inf_{i \leq k_n} \theta_0(i)}} + 1 - \frac{\inf_{h \in \mathcal{E}_{\theta_0, k_n}(M_n)} w_n\left(\theta_0 + \frac{h}{\sqrt{n}}\right)}{\sup_{h \in \mathcal{E}_{\theta_0, k_n}(M_n)} w_n\left(\theta_0 + \frac{h}{\sqrt{n}}\right)}\right).$$

If the sequence of truncation levels $(k_n)_n$ and radii $(M_n)_n$ satisfies Conditions (3.1) and (3.4) then

$$\mathbb{P}_{n,0}\left(\mathrm{NV}(M_n)\right) = o(1).$$

The third term $\left\|P_{H_n|X_{1:n}}^{M_n} - P_{H_n|X_{1:n}}\right\|$ is handled thanks to the posterior concentration lemma, since by Lemma B.1 in the appendix

$$\left\|P_{H_n|X_{1:n}}^{M_n} - P_{H_n|X_{1:n}}\right\| = 2 P_{H_n|X_{1:n}}\left(H_n^T I(\theta_0) H_n \geq M_n\right).$$

The proof of the Theorem is concluded by upper-bounding $\|\mathcal{N}_{k_n} - \mathcal{N}_{k_n}^{M_n}\|$. The latter quantity is a matter of concern because we are facing increasing dimensions $(k_n)_{n \in \mathbb{N}}$. It is checked in Section 4.4 that

**Proposition 3.11.** *There exists a universal constant $C$ such that if* $\liminf_n (n \inf_{i \leq k_n} \theta_0(i)) \geq c_0 > 0$ *and* $\liminf M_n/k_n \geq 64$, *then for large enough $n$*

$$\mathbb{P}_{n,0}\left\|\mathcal{N}_{k_n} - \mathcal{N}_{k_n}^{M_n}\right\| \leq C \exp\left(-\frac{M_n \wedge c_0 M_n^2}{C}\right)$$

□

### 3.2. Estimating functionals

The Bernstein-von Mises Theorem provides a handy tool to check the asymptotic normality of estimators of Rényi and Shannon entropies. Antos and Kontoyiannis (2001) established that plug-in estimators of Shannon and Rényi entropies are consistent whatever the sampling probability is. They also proved that entropy estimation may be arbitrarily slow, and that on a large class of sampling distributions, the mean squared error is $O(\log n/n)$. In the parametric setting, that is with fixed finite alphabets, analogues of the delta-method and the classical Bernstein-Von Mises Theorem can be used to check the asymptotic normality of both frequentist and Bayesian entropy estimators. Our purpose is to show that the non-parametric Bernstein-Von Mises Theorem can be used as well.

For any $\alpha > 0$, let $g_\alpha$ be the real function defined for non negative real numbers by $g_\alpha(u) = u^\alpha$ for $\alpha \neq 1$, and $g_1(u) = u \log u$ (with the convention $g_1(0) = 0$). The additive functional $G_\alpha$ is defined by

$$G_\alpha(\theta) = \sum_{i=1}^{+\infty} g_\alpha(\theta(i)).$$

The *Shannon entropy* of the probability mass function $\theta$ is $-G_1(\theta)$ and for $\alpha \neq 1$, $\frac{-1}{\alpha-1} \log G_\alpha(\theta)$ denotes the *Rényi entropy of order $\alpha$* (Cover and Thomas, 1991).

Let $T = (\tau(i))_{0 \leq i \leq k_n}$ be distributed according to the posterior distribution, a Bayesian estimator of $G_\alpha(\theta)$ may be constructed using the posterior distribution of

$$G_{n,\alpha}(T) = \sum_{i=1}^{k_n} g_\alpha(\tau(i)).$$

The Bernstein-Von Mises Theorem asserts that under $\mathbb{P}_{n,0}$, for large enough $n$, the posterior distribution of $(\tau(i))_{1 \leq i \leq k_n}$ is approximately Gaussian, centered around the maximum likelihood estimator $\widehat{\theta}_n = (\widehat{\theta}(i))_{1 \leq i \leq k_n}$, with variance $\frac{1}{n} I(\theta_0)^{-1}$. Theorem 3.12 below makes a similar assertion concerning $G_{n,\alpha}(T)$.

Let $G_{n,\alpha}(\widehat{\theta}_n)$ be the truncated plug-in maximum likelihood estimator:

$$G_{n,\alpha}(\widehat{\theta}_n) = \sum_{i=1}^{k_n} g_\alpha(\widehat{\theta}_n(i)).$$

The variance parameter $\gamma_{n,\alpha}$ is defined by

$$\gamma_{n,\alpha}^2 = \sum_{i=1}^{k_n} \theta_0(i) \left(g'_\alpha(\theta_0(i))\right)^2 - \left(\sum_{i=1}^{k_n} \theta_0(i) g'_\alpha(\theta_0(i))\right)^2.$$

Notice that

$$g'_\alpha(u) = \begin{cases} \alpha u^{\alpha-1} & \alpha \neq 1 \\ \log u + 1 & \alpha = 1 \end{cases}$$

so that $\gamma_{n,1}^2$ has limit $\gamma_1^2 = \sum_{i=1}^\infty \theta_0(i)(\log \theta_0(i) + 1)^2 - (\sum_{i=1}^\infty \theta_0(i)(\log \theta_0(i) + 1))^2$ as soon as this is finite, and $\gamma_{n,\alpha}^2$ has limit $\gamma_\alpha^2 = \alpha^2 [\sum_{i=1}^\infty \theta_0(i)^{2\alpha-1} - (\sum_{i=1}^\infty \theta_0(i)^\alpha)^2] = \alpha^2 [G_{2\alpha-1}(\theta_0) - (G_\alpha(\theta_0))^2]$ as soon as this is finite, which requires at least that $\alpha > \frac{1}{2}$.

Now, let $\mathcal{I}$ be the collection of all intervals in $\mathbb{R}$, and for any $I \in \mathcal{I}$, let $\Phi(I) = \int_I \phi(x)\mathrm{d}x$ where $\phi$ is the density of $\mathcal{N}(0,1)$. The following Theorem asserts that the Levy-Prokhorov distance between the posterior distribution of $\sqrt{n}\big(G_{n,\alpha}(T) - G_{n,\alpha}(\widehat{\theta}_n)\big)$ and $\mathcal{N}(0,\gamma_\alpha^2)$ tends to 0 in $\mathbb{P}_{n,0}$ probability. The Levy-Prokhorov distance metrizes convergence in distribution.

**Theorem 3.12.** (ESTIMATING FUNCTIONALS) *If $\lim_n \gamma_{n,\alpha}^2 = \gamma_\alpha^2$ is finite, then under the assumptions of the Bernstein-von Mises Theorem (Theorem 3.7),*

$$\sup_{I \in \mathcal{I}} \left| P_{H_n | X_{1:n}} \left( \frac{\sqrt{n}(G_{n,\alpha}(T) - G_{n,\alpha}(\widehat{\theta}_n))}{\gamma_{n,\alpha}} \in I \right) - \Phi(I) \right| \to 0$$

*in $\mathbb{P}_{n,0}$ probability.*

The proof of this theorem is given in Section 5.

Let us define the symmetric Bayesian credible set with would-be coverage probability $1 - \delta$ as the smallest interval which has posterior probability larger than $1 - \alpha$. This credible set is an empirical interval since it is defined thanks to an empirical quantity, the posterior distribution. In order to construct such a region, it is enough to sample from the posterior distribution using MCMC sampling methods. Note that this symmetric Bayesian credible set is not the (non fully empirical) interval

$$\left[ G_{n,\alpha}(\widehat{\theta}_n) - \frac{u_\delta \gamma_{n,\alpha}}{\sqrt{n}}; G_{n,\alpha}(\widehat{\theta}_n) + \frac{u_\delta \gamma_{n,\alpha}}{\sqrt{n}} \right]$$

where $u_\delta$ is the $1-\delta/2$ quantile of $\mathcal{N}(0,1)$. Theorem 3.12 just asserts that asymptotically, the symmetric Bayesian credible set has length $u_\delta \gamma_{n,\alpha}/\sqrt{n}$. and is centered around $G_{n,\alpha}(\widehat{\theta}_n)$. Hence Theorem 3.12 asserts that, in $\mathbb{P}_{n,0}$-probability, Bayesian credible sets for $G_\alpha(\theta_0)$ and frequentist confidence intervals based on truncated plug-in maximum likelihood estimators are asymptotically equivalent.

The next theorem provides sufficient conditions for the plug-in truncated maximum likelihood estimators to satisfy a central limit theorem with limiting variance $\gamma_\alpha^2$.

**Theorem 3.13.** (MLE FUNCTIONAL ESTIMATION) *Assume that $\lim_n \gamma_{n,\alpha}^2 = \gamma_\alpha^2$ is finite. If the truncation parameter $k_n$ satisfies:*

$$\begin{aligned}
(n \inf_{i \leq k_n} \theta_0(i))^{-1/2} k_n &= o(1) \\
\sum_{i=1}^{k_n} \theta_0(i) |g'_\alpha(\theta_0(i))|^3 &= o(\sqrt{n}) \\
\sum_{k_n+1}^{\infty} \theta_0(i) g_\alpha(\theta_0(i)) &= o\left(\frac{1}{\sqrt{n}}\right),
\end{aligned}$$

*then $\sqrt{n}\big(G_{n,\alpha}(\widehat{\theta}_n) - G_\alpha(\theta_0)\big)$ converges in distribution to $\mathcal{N}(0,\gamma_\alpha^2)$.*

### 3.3. Dirichlet prior distributions

We may now check that when using Dirichlet distributions as prior distributions, there exist truncation levels $(k_n)_n$ and radii $(M_n)_n$ such that Conditions 3.4 (prior smoothness) and 3.5 (prior concentration) hold.

Let $\boldsymbol{\beta} = (\beta_0, \beta_1, \ldots, \beta_{k_n})$ be a $(k_n + 1)$-tuple of positive real numbers. The Dirichlet distribution with parameter $(\beta_0, \beta_1, \ldots, \beta_{k_n})$ on the probability mass functions on $\{0, 1, \ldots, k_n\}$ has density

$$w_{n,\boldsymbol{\beta}}(\theta(1), \ldots, \theta(k_n)) = \frac{\Gamma\left(\sum_{i=0}^{k_n} \beta_i\right)}{\prod_{i=0}^{k_n} \Gamma(\beta_i)} \prod_{i=0}^{k_n} \theta(i)^{\beta_i - 1}.$$

In the absence of prior knowledge concerning the sampling distribution $\theta_0$, we refrain from assigning different masses on the coordinate components: we consider Dirichlet priors $W_{n,\beta}$ with constant parameter $\boldsymbol{\beta} = (\beta, \ldots, \beta)$ for some positive $\beta$.

Note that for $\beta = 1$ (the so-called Laplace prior), the Prior Smoothness Condition (3.4) trivially holds.

**Proposition 3.14.** *Let the sequence of prior distributions consist of the Dirichlet priors with parameter $\beta > 0$. The non-parametric Bernstein-Von Mises Theorem (3.7) holds if the sequence of truncation levels $(k_n)_{n \in \mathbb{N}}$, $k_n \to +\infty$, and the p.m.f. over $\mathbb{N}_*$, $\theta_0$ satisfy Condition 2.1, and if*

$$k_n \log n \vee \log\left(\det(I(\theta_0))\right) = o\left(\left(n \inf_{i \leq k_n} \theta_0(i)\right)^{1/3}\right).$$

Using such Dirichlet priors, checking the conditions of Theorem 3.7 boils down to checking the Prior Smoothness Condition.

*Proof.* For any $h$ and $g$ in $\mathcal{E}_{\theta_0, k_n}(M_n)$, for large enough $n$,

$$\sup_{h,g \in \mathcal{E}_{\theta_0,k_n}(M_n)} \frac{w_{n,\beta}\left(\theta_0 + \frac{h}{\sqrt{n}}\right)}{w_{n,\beta}\left(\theta_0 + \frac{g}{\sqrt{n}}\right)}$$

$$\leq \sup_{h,g \in \mathcal{E}_{\theta_0,k_n}(M_n)} \prod_{i=0}^{k_n} \left(\frac{\theta_0(i) + t\frac{h(i)}{\sqrt{n}}}{\theta_0(i) + \frac{g(i)}{\sqrt{n}}}\right)^{\beta-1}$$

$$\leq \sup_{h,g \in \mathcal{E}_{\theta_0,k_n}(M_n)} \exp\left[|\beta - 1| \sum_{i=0}^{k_n} \left\{\log\left(1 + \frac{|h(i)|}{\sqrt{n}\theta_0(i)}\right) - \log\left(1 - \frac{|g(i)|}{\sqrt{n}\theta_0(i)}\right)\right\}\right]$$

$$\leq \sup_{h,g \in \mathcal{E}_{\theta_0,k_n}(M_n)} \exp\left[|\beta - 1| \sum_{i=0}^{k_n} \left\{\frac{|h(i)|}{\sqrt{n}\theta_0(i)} + 2\frac{|g(i)|}{\sqrt{n}\theta_0(i)}\right\}\right]$$

$$\leq \exp 3\left[\sqrt{\frac{M_n(k_n + 1)(\beta - 1)^2}{n \inf_{i \leq k_n} \theta_0(i)}}\right],$$

as for each $g \in \mathcal{E}_{\theta_0,k_n}(M_n)$, $|g(i)|/\sqrt{n}\theta_0(i) \leq \sqrt{M_n/(n\inf_{i \leq k_n} \theta_0(i))}$ so that as soon as Condition 3.1 holds, for large enough $n$, $|g(i)|/\sqrt{n}\theta_0(i) \leq \frac{1}{2}$ and $\log(1 - |g(i)|/\sqrt{n}\theta_0(i)) \geq -2|g(i)|/\sqrt{n}\theta_0(i)$. Thus, the Prior Smoothness Condition holds as soon as

$$\frac{M_n k_n}{(n \inf_{i \leq k_n} \theta_0(i))} \to 0,$$

which is a consequence of Condition 3.1. On the other hand

$$-\log w_n(\theta_0) = O\left(\log(\det(I(\theta_0))) + k_n \log(k_n)\right).$$

Thus, using Dirichlet prior with parameter $\beta$, the Prior Smoothness and Prior Concentration Conditions hold for $\theta_0$ with truncation levels $k_n$ as soon as Condition 3.1 and

$$k_n \log n + \log(\det(I(\theta_0))) = o(M_n).$$

But the existence of a sequence of radii $(M_n)$ tending to infinity such that both the last condition and Condition 3.1 hold, is a straightforward consequence of Condition 2.1 and of the condition in Proposition 3.14. □

Note that if the prior distribution is Dirichlet with parameter $\boldsymbol{\beta}$ then the posterior distribution is Dirichlet with parameters $\boldsymbol{\beta} + (N_0, N_1, \ldots, N_{k_n})$. Let $n_i = \sum_{j<i} N_j$ for $i \leq k_n$, agreeing on $n_0 = 0$. Sampling from the posterior distribution is equivalent to picking an independent sample of $n$ exponentially distributed random variables, $Y_1, \ldots, Y_n$, picking another independent sample $Z_0, \ldots, Z_{k_n}$ of $k_n + 1$ independent $\Gamma(\beta, 1)$-distributed random variable, and letting $\theta^*(i) = \left(Z_i + \sum_{n_i < j \leq n_{i+1}} Y_j\right)/(\sum_{j=1}^n Y_j + \sum_{j=0}^{k_n} Z_j)$. The latter procedure is very close to the Bayesian Bootstrap (Rubin, 1981), indeed, we obtain the latter procedure if we omit to add the $Z_i$ in the weights. This procedure which has been extensively investigated (See Lo, 1988, 1987; Weng, 1989, among other references) is now considered as a special case of exchangeable bootstrap (See van der Vaart and Wellner, 1996, and references therein). Theorems from the preceding section tell us that the Bayesian bootstrap of (non-linear) functionals of the sampling distribution approximate the asymptotic distribution of maximum likelihood estimates. We leave the analysis of the second-order properties of the posterior distribution to further investigations.

### 3.4. Examples

Previous results may now be applied to two examples of envelope classes already investigated by Boucheron et al. (2009):

1. The sampling probability $\theta_0$ is said to have exponential($\eta$) decay if there exists $\eta > 0$, and a positive constant $C$ such that

$$\forall i \in \mathbb{N}_*, \ \frac{1}{C}\exp(-\eta i) \leq \theta_0(i) \leq C\exp(-\eta i).$$

Using truncation level $k_n$,

$$\frac{\exp(-\eta)}{C(1-\exp(-\eta))}\exp(-\eta k_n) \leq \theta_0(0) \leq \frac{C\exp(-\eta)}{1-\exp(-\eta)}\exp(-\eta k_n).$$

2. The sampling probability $\theta_0$ is said to have polynomial($\eta$) decay if there exists $\eta > 1$, and a positive constant $C$ such that

$$\forall i \in \mathbb{N}_*, \; \frac{1}{C\,i^\eta} \leq \theta_0(i) \leq \frac{C}{i^\eta}.$$

Using truncation level $k_n$, $\frac{c}{(k_n+1)^{\eta-1}} \leq \theta_0(0) \leq \frac{C}{k_n^{\eta-1}}$.

Let us first assume that $\theta_0$ has exponential($\eta$) decay. Then with $\tilde{c} = \frac{1}{C} \wedge \frac{\exp(-\eta)}{C(1-\exp(-\eta))}$,

$$\inf_{i\leq k_n} \theta_0(i) \geq \tilde{c}\exp(-\eta k_n),$$

$$-\sum_{i=0}^{k_n} \log\theta_0(i) \leq \left(\eta\frac{k_n(k_n+3)}{2} - (k_n+1)\log C - \log\frac{\exp(-\eta)}{1-\exp(-\eta)}\right).$$

Invoking Proposition 3.14, the non-parametric Bernstein-Von Mises Theorem holds for $\theta_0$ with exponential($\eta$) decay using the Dirichlet prior with parameter $\beta > 0$ with truncation levels

$$k_n = \frac{1}{\eta}(\log n - a\log\log n), \quad a > 6.$$

Theorems 3.12 and 3.13 apply as soon as $\alpha > \frac{1}{2}$, so that the Bayesian estimates of entropy and of Rényi-entropy of order $\alpha > \frac{1}{2}$ satisfy a Bernstein-von-Mises theorem with $\sqrt{n}$-rate.

When $\theta_0$ has polynomial($\eta$) decay, $\inf_{i\leq k_n}\theta_0(i) \geq 1/(C\,k_n^\eta)$, and

$$-\sum_{i=0}^{k_n}\log\theta_0(i) \leq \left(\eta k_n\log k_n - (k_n+1)\log C\right).$$

Invoking Proposition 3.14, the non-parametric Bernstein-Von Mises Theorem holds using the Dirichlet prior with parameter $\beta > 0$ with truncation levels

$$k_n = \left(\frac{n}{u_n}\right)^{\frac{1}{\eta+3}}, \quad \frac{u_n}{(\log n)^3} \to +\infty.$$

Theorems 3.12 and 3.13 concerning estimations of functionals hold as soon as $2\alpha > 1 + 1/\eta$, so that the Bayesian estimates of entropy and of Rényi-entropy of order $\alpha > 1/2 + 1/(2\eta)$ satisfy a Bernstein-von-Mises theorem with rate $\sqrt{n}$.

### 3.5. Comparison with Ghosal's conditions

Now, we aim at comparing the set of conditions used by Ghosal (2000) to establish a Bernstein-Von Mises Theorem for sequences of multinomial models using log-odds parametrization. An exhaustive comparison of the two approaches (that is, comparing the merits of combining Le Cam's proof and concentration inequalities for some quadratic forms with the merits of Ghosal's proof which refines Portnoy's arguments) should first be based on a general purpose result characterizing the impact of re-parametrization on asymptotic normality of posterior distributions. This would exceed the ambitions of this paper. Then a thorough comparison between conditions (P) (Prior Smoothness and Concentration) and (R) (Prior concentration and behavior of likelihood ratios in the vicinity of the target $\theta_0$) and the conditions used in this paper would be in order. As a matter of fact, provided re-parametrization is taken into account, the prior smoothness conditions in the two papers are not essentially different. On the other hand the conditions on the integrability of likelihood ratios seem somewhat different. Looking for general exponential families, Ghosal (2000) imposes upper-bounds on the fourth and the third moment of linear forms of $\sqrt{I(\theta)}\Delta_1(\theta)$ for $\theta$ close to $\theta_0$ (this is the meaning of conditions on the growth of $B_{1,n}(c)$ and $B_{2,n}(c)$.) In this paper, we take advantage of the fact that $\Delta_n(\theta)$ is a multinomial vector.

Keep in mind that we refrain from assuming that all $\theta_0(i), i \leq k_n$ are of order $1/k_n$ as in (Ghosal, 2000, page 60). Indeed, we consider situations where $k_n \inf_{i \leq k_n} \theta_0(i) = o(1)$ as in Section 3.4. The trace of the information matrix $I(\theta_0)$ (which coincides with $\mathbf{F}^{-1}$ using Ghosal's notations) is equal to $\sum_{i=1}^{k_n} 1/\theta_0(i) + k_n/\theta_0(0) \leq 2k_n/\inf_{i \leq k_n} \theta_0(i)$, and it may not be $O(k_n^2)$, as in (Ghosal, 2000, page 60). For example, using the notations from Section 3.4, if $\theta_0$ has polynomial-($\eta$) decay, $\sum_{i=1}^{k_n} 1/\theta_0(i) + k_n/\theta_0(0) \geq \frac{1}{C} \sum_{i=1}^{k_n} i^\eta + \frac{1}{C} k_n^{\eta-1} \geq \frac{ck_n^{\eta+1}}{\eta+1} + \frac{k_n^{\eta-1}}{C}$.

In this setting, we may even look at the growth of $B_{2,n}(0)$ (as defined in (Ghosal, 2000)) as $n$ tends to infinity

$$B_{2,n}(c) = \sup \left\{ \mathbb{P}_\theta \left[ \left| \mathbf{a}^T I^{\frac{1}{2}}(\theta) \Delta_n(\theta) \right|^4 \right] ; \|\mathbf{a}\| = 1, \ \operatorname{Var}_{\theta_0} \left( \log \frac{\theta(i)}{\theta_0(i)} \right) \leq \frac{ck_n}{n} \right\}.$$

Choosing $\mathbf{a}$ as $\frac{1}{\sqrt{k_n}} \mathbf{1}$ and carefully performing straightforward computations, it is possible to check that if $\theta_0$ has polynomial-($\eta$) decay (according to the framework of Section 3.4), $B_{2,n}(0) \geq C k_n^{2\eta}$, so that the clause $B_{2,n}(c \log k_n) k_n^2 (\log k_n)/n \to 0$ for all $c > 0$ in Condition (R), implies $k_n^{2+2\eta} \log k_n/n \to 0$. This condition is more demanding that the conditions we obtained at the end of Section 3.4.

### 3.6. Classical non-parametric approach to posterior concentration

We compare the posterior concentration lemma (Lemma 3.6) and the classical results on posterior concentration obtained in non-parametric statistics (See

Ghosal and van der Vaart (2007a),Ghosal et al. (2000), Ghosal and van der Vaart (2007b),Ghosal and van der Vaart (2001)).

Let $\Theta_{k_n}$ denote the set of probability distributions over $\{0, \ldots, k_n\}$. Let $\epsilon_n^2$ satisfy $n\epsilon^2 = M_n$.

Let $V_n(\epsilon_n)$ be the set:

$$V_n(\epsilon_n) = \left\{ \theta \ : \ \sum_{i=0}^{k_n} \theta_0(i) \log \frac{\theta_0(i)}{\theta(i)} \leq \epsilon_n^2 \text{ and } \sum_{i=0}^{k_n} \theta_0(i) \left( \log \frac{\theta_0(i)}{\theta(i)} \right)^2 \leq \epsilon_n^2 \right\}.$$

Let $d$ denote the Hellinger distance between probability mass functions:

$$d(\theta_1, \theta_2) = \left[ \sum_{i=0}^{k_n} \left( \sqrt{\theta_1(i)} - \sqrt{\theta_2(i)} \right)^2 \right]^{1/2}$$

Let $D(\epsilon, \Theta_{k_n}, d)$ denote the $\epsilon$-packing number of $\Theta_{k_n}$, that is the maximum number of points in $\Theta_{k_n}$ such that the Hellinger distance between every pair is at least $\epsilon$.

Theorem 2.1 in Ghosal et al. (2000) asserts that, if for some $C > 0$, we have $W_n\{V_n(\epsilon_n)\} \geq \exp(-Cn\epsilon_n^2)$ and if $\log D(\epsilon_n, \Theta_{k_n}, d) \leq n\epsilon_n^2$, then for large enough $A$,

$$\mathbb{P}_{\cdot | X_{1:n}}\{d(\theta, \theta_0) \geq A\epsilon_n\} \to 0$$

in $\mathbb{P}_{0,n}$ probability.

In this paper, the prior $W_n$ is supported by $\Theta_{k_n}$, and a careful reading shows that the proof in Ghosal et al. (2000) can be adapted to situations where the sampling probability changes with $n$.

Now, $\Theta_{k_n}$ endowed with the Hellinger distance is isometric to the intersection of the positive quadrant and the unit ball of $\mathbb{R}^{k_n+1}$ endowed with the Euclidean metric, so that there exists a universal constant $C$

$$\frac{1}{C} \cdot \left( \frac{1}{2\epsilon} \right)^{k_n} \leq D(\epsilon, \Theta_{k_n}, d) \leq C \cdot \left( \frac{1}{2\epsilon} \right)^{k_n}$$

and $\log D(\epsilon_n, \Theta_{k_n}, d) \leq n\epsilon_n^2$ if and only if $k_n \log \frac{n}{M_n} \leq CM_n$.

If $h^T I(\theta_0) h \leq M_n$, then $\sup_{i \leq k_n} \left| \frac{h(i)}{\sqrt{n}\theta_0(i)} \right| \leq \sqrt{\frac{M_n}{n \inf_{i \leq k_n} \theta_0(i)}}$. Hence, for $i \leq k_n$

$$\log \left( 1 + \frac{h(i)}{\sqrt{n}\theta_0(i)} \right) = \frac{h(i)}{\sqrt{n}\theta_0(i)} - \frac{h(i)^2}{2n\theta_0(i)^2} + o\left( \frac{h(i)^2}{n\theta_0(i)^2} \right)$$

So that, letting $\theta(i) = \theta_0(i) + \frac{h(i)}{\sqrt{n}}$,

$$\sum_{i=0}^{k_n} \theta_0(i) \log \frac{\theta_0(i)}{\theta(i)} = \frac{\sigma_n^2(h)}{2n} + o\left( \frac{\sigma_n^2(h)}{n} \right)$$

and

$$\sum_{i=0}^{k_n} \theta_0(i) \left( \log \frac{\theta_0(i)}{\theta(i)} \right)^2 = \frac{\sigma_n^2(h)}{n} + o\left( \frac{\sigma_n^2(h)}{n} \right)$$

Hence, as soon as $\frac{M_n}{n \inf_{i \leq k_n} \theta_0(i)} = o(1)$, if for some $\delta > 0$ and some $c > 0$,

$$W_n \left\{ h \ : \ h^T I(\theta_0) h \leq \delta M_n \right\} \geq e^{-cM_n}$$

then for some constant $C > 0$,

$$W_n \left\{ V_n \left( \epsilon_n \right) \right\} \geq e^{-CM_n}.$$

Under our assumptions, the non-parametric Bayesian theorem implies that for large enough $A$, $\mathbb{P}_{\cdot | X_{1:n}} \left\{ d\left(\theta, \theta_0\right) \geq A\epsilon_n \right\} \to 0$ in $\mathbb{P}_{0,n}$ probability.

However, Lemma 3.6 posterior concentration with respect to the Fisher distance:

$$\mathbb{P}_{\cdot | X_{1:n}} \left\{ (\theta - \theta_0)^T I(\theta_0) (\theta - \theta_0) \geq \epsilon_n^2 \right\} \to 0$$

in $\mathbb{P}_{0,n}$ probability.

As the Fisher distance upper-bounds the squared Hellinger distance (See Tsybakov, 2004), Lemma 3.6 implies the generic posterior concentration lemma. But Lemma 3.6 could not be deduced from generic posterior concentration lemma since the Fisher distance cannot be upper-bounded by a linear function of the Hellinger distance. As a matter of fact,

$$(\theta - \theta_0)^T I(\theta_0) (\theta - \theta_0) \ \leq \ d(\theta, \theta_0) \left( 1 + \frac{1}{\inf_{0 \leq i \leq k_n} \theta_0(i)} \right).$$

Hence, if $d(\theta, \theta_0) = o(\inf_{0 \leq i \leq k_n} \theta_0(i))$, Hellinger and Fisher metrics are comparable, but this does not hold in full generality. For instance, for $\theta$ such that $\theta(k_n) = \theta_0(k_n) + \sqrt{\theta_0(k_n)}$, $\theta(1) = \theta_0(k_n) - \sqrt{\theta_0(k_n)}$, and $\theta(i) = 0$ for $i \neq 1$ and $i \neq k_n$, then $d(\theta, \theta_0) \to 0$, but $(\theta - \theta_0)^T I(\theta_0) (\theta - \theta_0) \sim 1$.

## 4. Proof of the Bernstein-Von Mises Theorem

In this section, we establish the building blocks of the proof of the Bernstein-Von Mises Theorem that is Proposition 3.10, the posterior concentration Lemma and Proposition 3.11.

### 4.1. Truncated distributions

In order to prove Proposition 3.10, it is enough to upper bound $\mathrm{NV}(M_n)$:

$$\iint \left( 1 - \frac{w_n(\theta_0 + \frac{g}{\sqrt{n}})}{w_n(\theta_0 + \frac{h}{\sqrt{n}})} e^{\left\{ \frac{A_n(g) - A_n(h)}{2} + C_n(g) - C_n(h) \right\}} \right)^+ \mathrm{d}\mathcal{N}_{k_n}^{M_n}(g) \mathrm{d}P_{H_n | X_{1:n}}^{M_n}(h),$$

where $A_n$ and $C_n$ are defined in Section 3.1.

We take advantage of the fact that integration is performed on $\mathcal{E}_{\theta_0, k_n}(M_n)$, in order to uniformly upper-bound the integrand.

Using the duality between $\ell^1_{k_n}$ and $\ell^\infty_{k_n}$, for all $h \in \mathcal{E}_{\theta_0,k_n}(M_n)$

$$A_n(h) \leq \sup_{h \in \ell^1_{k_n}, \|h\|_1 \leq M_n} \sum_{i=0}^{k_n} h(i) \left( \frac{N_i}{n\theta_0(i)} - 1 \right) = M_n \sup_{i=0,\ldots,k_n} \left| \frac{N_i}{n\theta_0(i)} - 1 \right|$$

and

$$|C_n(h)| \leq M_n \sup_{i=0,\ldots,k_n} \left| \frac{N_i}{n\theta_0(i)} \right| \left| R \left( \frac{\sqrt{M_n}}{\sqrt{n \inf_{i \leq k_n} \theta_0(i)}} \right) \right|.$$

Then as $(1 - (1-x)e^{-y})^+ \leq x + y$ for $x, y \geq 0$,

$$\mathrm{NV}(M_n) \leq M_n \times \left( \sup_{i \leq k_n} \left| \frac{N_i}{n\theta_0(i)} - 1 \right| + 2 \sup_{i \leq k_n} \left| \frac{N_i}{n\theta_0(i)} \right| R \left( \frac{\sqrt{M_n}}{\sqrt{n \inf_{i \leq k_n} \theta_0(i)}} \right) \right)$$
$$+ \left( 1 - \frac{\inf_{h \in \mathcal{E}_{\theta_0,k_n}(M_n)} w_n \left( \theta_0 + \frac{h}{\sqrt{n}} \right)}{\sup_{h \in \mathcal{E}_{\theta_0,k_n}(M_n)} w_n \left( \theta_0 + \frac{h}{\sqrt{n}} \right)} \right).$$

The second term can be upper-bounded assuming the prior smoothness condition. The first term is a sum of two random suprema.

The expected value of the maximum of random variables with uniformly controlled logarithmic moment generating functions can be handily upper-bounded thanks to an argument due to Pisier (Massart, 2003): if $(W_i)_{1 \leq i \leq k}$ are real random variables, then

$$\mathbb{P}\left[ \sup_{i \leq k} W_i \right] \leq \inf_{\lambda > 0} \frac{1}{\lambda} \left\{ \log k + \sup_{i \leq k} \log \mathbb{P}\left[ \exp \lambda W_i \right] \right\}. \tag{4.1}$$

For each $i$, the random variable $N_i$ is binomially distributed with parameters $n$ and $\theta_0(i)$, $\log(1+u) \leq u$, and for all $u \geq 0$, $\frac{e^u - 1}{u} - 1 \geq \frac{e^{-u} - 1}{u} + 1$, so that (4.1) leads to

$$\mathbb{P}_{\theta_0}\left[ \sup_{i \leq k_n} \left| \frac{N_i}{n\theta_0(i)} - 1 \right| \right] \leq \inf_{\lambda > 0} \left\{ \frac{\log 2(k_n + 1)}{\lambda} - 1 + \sup_{i \leq k_n} \frac{\exp\left(\frac{\lambda}{n\theta_0(i)}\right) - 1}{\frac{\lambda}{n\theta_0(i)}} \right\},$$

so that choosing $\lambda = \sqrt{\log(2(k_n+1)) n \inf_{i \leq k_n} \theta_0(i)}$, as the function $u \to \frac{e^u - 1}{u} - 1$ is increasing on $\mathbb{R}_+$, letting $\delta_n = \sqrt{\frac{\log(2(k_n+1))}{n \inf_{i \leq k_n} \theta_0(i)}}$,

$$\mathbb{P}_{\theta_0}\left[ \sup_{i \leq k_n} \frac{N_i}{n\theta_0(i)} \right] - 1 \leq \mathbb{P}_{\theta_0}\left[ \sup_{i \leq k_n} \left| \frac{N_i}{n\theta_0(i)} - 1 \right| \right] \leq \delta_n + \frac{\exp(\delta_n) - 1}{\delta_n} - 1.$$

Thus,

$$\mathbb{P}_{n,0}\mathrm{NV}(M_n) \leq M_n \left( \left( \delta_n + \frac{\exp(\delta_n) - 1}{\delta_n} \right) \left( 1 + R\left( \frac{\sqrt{M_n}}{\sqrt{n \inf_{i \leq k_n} \theta_0(i)}} \right) \right) - 1 \right)$$
$$+ 1 - \frac{\inf_{h \in \mathcal{E}_{\theta_0,k_n}(M_n)} w_n \left( \theta_0 + \frac{h}{\sqrt{n}} \right)}{\sup_{h \in \mathcal{E}_{\theta_0,k_n}(M_n)} w_n \left( \theta_0 + \frac{h}{\sqrt{n}} \right)}$$

and the proposition follows using Assumptions (3.1) and (3.4) and the fact that $R(u) = O(u)$ as $u$ tends toward 0. □

### 4.2. Tail bounds for quadratic forms

In this section, we gather a few results concerning tail bounds for quadratic forms or square-roots of quadratic forms in Gaussian and empirical settings. All those bounds are obtained by resorting to concentration inequalities for Gaussian distributions or for suprema of empirical processes.

Let us first start by a first bound concerning chi-square distributions. Let $\xi_n^2$ be distributed according to $\chi_{k_n}^2$ (chi-square distribution with $k_n$ degrees of freedom), the following inequality is a direct consequence of Cirelson's inequality (Massart, 2003):

$$\mathbb{P}\left\{\xi_n \geq \sqrt{k_n} + \sqrt{2x}\right\} \leq \exp(-x). \tag{4.2}$$

The following handy inequality provides non-asymptotic tail-bounds for Pearson statistics. For any $\theta \in \Theta_{k_n}$ let $V_n(\theta)$ denote the square root of the Pearson statistic

$$V_n(\theta) = \left(\sum_{i=0}^{k_n} \frac{(N_i - n\theta(i))^2}{n\theta(i)}\right)^{1/2} = \left(\Delta_n^T(\theta) I(\theta) \Delta_n(\theta)\right)^{1/2}.$$

The following follows from Talagrand's inequality for suprema of empirical processes (Massart, 2003, p. 170): for all $x > 0$,

$$\mathbb{P}_{n,0}\left(V_n(\theta_0) \geq 2\sqrt{k_n} + \sqrt{2x} + 3\frac{x}{\sqrt{n \inf_{i \leq k_n} \theta_0(i)}}\right) \leq \exp(-x). \tag{4.3}$$

Non-centered Pearson statistics also show up while proving the posterior concentration lemma. Let $\theta = \theta_0 + \frac{h}{\sqrt{n}}$ with $\sigma_n(h) \geq \sqrt{M_n}$. Note that from the definition of $V_n(\theta_0)$, it follows that

$$\begin{aligned}
V_n(\theta_0) &= \sup_{\mathbf{a}:\|\mathbf{a}\|=1} \sum_{i=0}^{k_n} \mathbf{a}_i \frac{N_i - n\theta_0(i)}{\sqrt{n\theta_0(i)}} \\
&= \sup_{\mathbf{a}:\|\mathbf{a}\|=1} \sum_{i=0}^{k_n} \mathbf{a}_i \left(\frac{N_i - n\theta(i)}{\sqrt{n\theta_0(i)}} + \sqrt{n}\frac{\theta(i) - \theta_0(i)}{\sqrt{\theta_0(i)}}\right) \\
&\geq \sum_{i=0}^{k_n} \mathbf{a}_i^* \frac{N_i - n\theta(i)}{\sqrt{n\theta_0(i)}} + \sigma_n(h),
\end{aligned}$$

where $\mathbf{a}_i^* = \frac{h(i)}{\sqrt{\theta_0(i)\sigma^2(h)}}$ for all $i \leq k_n$. So that

$$\mathbb{P}_{n,h}\left(V_n(\theta_0) \leq s_n\right) \leq \mathbb{P}_{n,h}\left(\sum_{i=0}^{k_n} \mathbf{a}_i^* \frac{N_i - n\theta(i)}{\sqrt{n\theta_0(i)}} \leq -\sigma_n(h) + s_n\right).$$

Computations carried out in the Appendix allow to establish that if $\sigma_n^2(h) \geq M_n$, and if $M_n = o(n \inf_{i \leq k_n} \theta_0(i))$,

$$\mathbb{P}_{n,h}\left\{V_n(\theta_0) < 2\sqrt{k_n} + \sqrt{\frac{M_n}{2}} + \frac{3M_n}{4\sqrt{n\inf_{i\leq k_n}\theta_0(i)}}\right\} \leq 2\exp\left(-\frac{M_n}{96}\right). \quad (4.4)$$

### 4.3. Proof of the posterior concentration lemma

*Proof.* We need to check that $P_{H_n|X_{1:n}}\left\{H_n^T I(\theta_0) H_n \geq M_n\right\}$ is small in $\mathbb{P}_{n,0}$ probability. For any $\theta \in \Theta_{k_n}$, let $V_n(\theta)$ denote the square root of the Pearson statistic

$$V_n(\theta) = \left(\sum_{i=0}^{k_n} \frac{(N_i - n\theta(i))^2}{n\theta(i)}\right)^{1/2} = \left(\Delta_n^T(\theta) I(\theta) \Delta_n(\theta)\right)^{1/2}.$$

A sequence of tests $(\phi_n)_{n \in \mathbb{N}}$ is defined by $\phi_n = \mathbb{1}_{V_n(\theta_0) \geq s_n}$ where each threshold $s_n$ is defined by $s_n = 2\sqrt{k_n} + \sqrt{2x_n} + 3x_n/\sqrt{n\inf_{i\leq k_n}\theta_0(i)}$ with $x_n = \frac{M_n}{4}$. The tests $\phi_n$ aim at separating $\theta_0$ from the complements of Fisher balls centered at $\theta_0$, that is from $\left\{\theta_0 + h\sqrt{n} : \sigma_n^2(h) \geq M_n\right\}$.

Hence, we need to check that

$$\begin{aligned}
&P_{H_n|X_{1:n}}\left\{H_n^T I(\theta_0) H_n \geq M_n\right\} \\
&= P_{H_n|X_{1:n}}\left\{H_n^T I(\theta_0) H_n \geq M_n\right\} \phi_n \\
&\quad + P_{H_n|X_{1:n}}\left\{H_n^T I(\theta_0) H_n \geq M_n\right\} (1 - \phi_n) \\
&\leq \phi_n + P_{H_n|X_{1:n}}\left\{H_n^T I(\theta_0) H_n \geq M_n\right\} (1 - \phi_n).
\end{aligned}$$

is small in $\mathbb{P}_{n,0}$ probability. As, in order to upper-bound $\mathbb{P}_{n,0}\phi_n$, it is enough to bound the tail of Pearson's statistics under $\mathbb{P}_{n,0}$, we focus on the expected value of the second term. Note that the latter is null as soon as the maximum likelihood estimator errs too far away from $\theta_0$.

In order to control $\mathbb{P}_{n,0} P_{H_n|X_{1:n}} \left\{H_n^T I(\theta_0) H_n \geq M_n\right\} (1 - \phi_n)$, we resort to the same contiguity trick as in (van der Vaart, 1998). Let $A$ be a fixed positive real, define the probability distribution $\mathbb{P}_{n,A}$ on $\mathbb{N}^n$ as the mixture of $\mathbb{P}_{n,h}$ when the prior is conditioned on the ellipsoid $\theta_0 + \frac{1}{\sqrt{n}}\mathcal{E}_{\theta_0,k_n}(A)$:

$$\mathbb{P}_{n,A}(B) = \frac{\int_{\mathcal{E}_{\theta_0,k_n}(A)} \mathbb{P}_{n,h}(B) w_n\left(\theta_0 + \frac{h}{\sqrt{n}}\right) \mathrm{d}h}{W_n\left(\theta_0 + \frac{1}{\sqrt{n}}\mathcal{E}_{\theta_0,k_n}(A)\right)}.$$

Arguing as in (van der Vaart, 2002), thanks to Lemma 2.3, one can check that the sequences $(\mathbb{P}_{n,0})_n$ and $(\mathbb{P}_{n,A})_n$ are mutually contiguous (for the sake of self-reference, a proof is given in the Appendix, see Section A). Hence, it is enough

to upper-bound

$$\mathbb{P}_{n,A}\left[P_{H_n|X_{1:n}}\left\{H_n^T I(\theta_0)H_n \geq M_n\right\}(1-\phi_n)\right]$$
$$\leq \frac{1}{W_n\{\theta_0 + \mathcal{E}_{0,k_n}(A)/\sqrt{n}\}} \int_{h \notin \mathcal{E}_{0,k_n}(M_n)} \mathbb{P}_{n,h}(1-\phi_n) w_n\left(\theta_0 + \frac{h}{\sqrt{n}}\right) \mathrm{d}h$$
$$\leq \frac{\sup_{h^T I(\theta_0)h \geq M_n} \mathbb{P}_{n,h}(1-\phi_n)}{W_n\{\theta_0 + \mathcal{E}_{\theta_0,k_n}(A)/\sqrt{n}\}}.$$

We will handle $\mathbb{P}_{n,0}\phi_n$ and $\mathbb{P}_{n,h}(1-\phi_n)$ using non-asymptotic upper bounds for centered and non-centered Pearson statistics while the prior mass around $\theta_0$ ($W_n\{\theta_0 + \mathcal{E}_{\theta_0,k_n}(A)/\sqrt{n}\}$) can be lower-bounded by assuming Conditions 3.4 and 3.5.

A direct application of Inequality (4.3) gives $\mathbb{P}_{n,0}\phi_n = \mathbb{P}_{n,0}\{V_n(\theta_0) \geq s_n\} \leq \exp(-x_n) = \exp(-M_n/4)$.

Non-centered Pearson statistics show up while handling $\mathbb{P}_{n,h}(1-\phi_n)$. Indeed $\mathbb{P}_{n,h}(1-\phi_n) = \mathbb{P}_{n,h}(V_n(\theta_0) \leq s_n)$. Let $\theta = \theta_0 + \frac{h}{\sqrt{n}}$ with $\sigma_n(h) \geq \sqrt{M_n}$. Then, using the definition of $\phi_n$, Inequality (4.4) entails

$$\mathbb{P}_{n,h}(1-\phi_n) \leq 2\exp\left(-\frac{M_n}{96}\right).$$

Let us now lower bound $W_n\{\theta_0 + \mathcal{E}_{\theta_0,k_n}(A)/\sqrt{n}\}$. Performing a change of variables (agreeing on the convention that $h(0) = -\sum_{i=1}^{k_n} h(i)$) leads to

$$W_{n,\alpha}\{n(\theta-\theta_0)^T I(\theta_0)(\theta-\theta_0) \leq A\}$$
$$= \int_{h \in \mathcal{E}_{\theta_0,k_n}(A)} w_n\left(\theta_0 + \frac{h}{\sqrt{n}}\right)\left(\frac{1}{\sqrt{n}}\right)^{k_n} \prod_{i=1}^{k_n} \mathrm{d}h(i)$$
$$\geq \left(\frac{\sup_{h \in \mathcal{E}_{\theta_0,k_n}(A)} w_n\left(\theta_0 + \frac{h}{\sqrt{n}}\right)}{\inf_{h \in \mathcal{E}_{\theta_0,k_n}(A)} w_n\left(\theta_0 + \frac{h}{\sqrt{n}}\right)}\right)^{-1} w_n(\theta_0) \left(\frac{1}{n}\right)^{\frac{k_n}{2}} \int_{h \in \mathcal{E}_{\theta_0,k_n}(A)} \prod_{i=1}^{k_n} \mathrm{d}h(i).$$

But the volume of the ellipsoid in $\mathbb{R}^{k_n}$ induced by $\mathcal{E}_{\theta_0,k_n}(A)$ is the inverse of the square root of the determinant of $I(\theta_0)$ (that is $\prod_{i=0}^{k_n} \theta_0(i)^{1/2}$) times the volume of the sphere with radius $\sqrt{A}$ in $\mathbb{R}^{k_n}$, that is $A^{k_n/2} \frac{2\Gamma(\frac{1}{2})^{k_n}}{k_n \Gamma(\frac{k_n}{2})}$ so that

$$W_{n,\alpha}\{n(\theta-\theta_0)^T I(\theta_0)(\theta-\theta_0) \leq A\}$$
$$\geq \left(\frac{\sup_{h \in \mathcal{E}_{\theta_0,k_n}(A)} w_n\left(\theta_0 + \frac{h}{\sqrt{n}}\right)}{\inf_{h \in \mathcal{E}_{\theta_0,k_n}(A)} w_n\left(\theta_0 + \frac{h}{\sqrt{n}}\right)}\right)^{-1} w_n(\theta_0) \prod_{i=0}^{k_n} \theta_0(i)^{1/2} \left(\frac{A}{n}\right)^{\frac{k_n}{2}} \frac{2\Gamma(\frac{1}{2})^{k_n}}{k_n \Gamma(\frac{k_n}{2})}.$$

Thus, assuming conditions 3.4 and 3.5:

$$\frac{\sup_{h^T I(\theta_0)h \geq M_n} \mathbb{P}_{n,h}(1-\phi_n)}{W_n\{\theta_0 + \mathcal{E}_{\theta_0,k_n}(A)/\sqrt{n}\}} \leq C\exp\left(-\frac{M_n}{96}\right)(1+o(1)).$$

□

### 4.4. Posterior Gaussian concentration

Proving Proposition 3.11 amounts to checking that the growth rate of the sequence of radii $M_n$ is large enough so as to balance the growth rate of dimension $k_n$. By Lemma B.1:

$$\left\|\mathcal{N}_{k_n}^{M_n} - \mathcal{N}_{k_n}\right\| \;=\; 2\int_{\sigma_n(h) \geq M_n} \mathrm{d}\mathcal{N}_{k_n}(\Delta_{n,\theta_0}, I^{-1}(\theta_0))(h).$$

The right-hand-side can be upper-bounded:

$$\begin{aligned}
&\int_{\sigma_n(h) \geq M_n} \mathrm{d}\mathcal{N}_{k_n}(\Delta_{n,\theta_0}, I^{-1}(\theta_0))(h) \\
&= \int \mathbf{1}_{(h+\Delta_{n,\theta_0})^T I(\theta_0)(h+\Delta_{n,\theta_0}) \geq M_n} \mathrm{d}\mathcal{N}_{k_n}(0, I^{-1}(\theta_0))(h) \\
&\leq \int \mathbf{1}_{2h^T I(\theta_0)h + 2\Delta_{n,\theta_0}^T I(\theta_0)\Delta_{n,\theta_0} \geq M_n} \mathrm{d}\mathcal{N}_{k_n}(0, I^{-1}(\theta_0))(h) \\
&\leq \int \mathbf{1}_{h^T I(\theta_0)h \geq M_n/4} \mathrm{d}\mathcal{N}_{k_n}(0, I^{-1}(\theta_0))(h) + \mathbf{1}_{\Delta_{n,\theta_0}^T I(\theta_0)\Delta_{n,\theta_0} \geq M_n/4} \\
&= \int \mathbf{1}_{\|h\|_2 \geq \sqrt{M_n}/2} \mathrm{d}\mathcal{N}_{k_n}(0, \mathrm{Id}_{k_n})(h) + \mathbf{1}_{\Delta_{n,\theta_0}^T I(\theta_0)\Delta_{n,\theta_0} \geq M_n/4}.
\end{aligned}$$

so that

$$\mathbb{P}_{n,0}\left\|\mathcal{N}_{k_n}^{M_n} - \mathcal{N}_{k_n}\right\| \;\leq\; P\left\{\xi_n \geq \sqrt{\tfrac{M_n}{4}}\right\} + \mathbb{P}_{n,0}\left(\Delta_{n,\theta_0}^T I(\theta_0)\Delta_{n,\theta_0} \geq \tfrac{M_n}{4}\right)$$

where $\xi_n^2$ is distributed according to $\chi_{k_n}^2$ (chi-square distribution with $k_n$ degrees of freedom). Then, invoking (4.2),

$$P\left\{\xi_n \geq \sqrt{\tfrac{M_n}{4}}\right\} \leq \exp\left(-\tfrac{M_n}{32}\right).$$

The second term in the upper bound is handled using (4.3) and choosing $x = \inf\left(\tfrac{M_n}{128}, \tfrac{c_0 M_n^2}{512}\right)$.

### 5. Proof of Theorem 3.12

In frequentist statistics, once asymptotic normality has been proved for an estimator, the so-called delta-method allows to extend this result to smooth functionals of this estimator. In this section, we develop an *ad hoc* approach that parallels the classical derivation of the delta-method. Taylor expansions allow to write $\sqrt{n}(G_{n,\alpha}(T) - G_\alpha(\theta_0))$ as the sum of a linear function of $H_n - \Delta_n(\theta_0)$ and of two (random) quadratic forms. Checking the theorem amounts to establish that under $\mathbb{P}_{n,0}$ those two quadratic forms converge to 0 in distribution.

Recall that $H_n = \sqrt{n}(\tau(i) - \theta_0(i))_{i=1}^{k_n}$ and $\Delta_n(\theta_0) = \sqrt{n}(\hat{\theta}(i) - \theta_0(i))_{i=1}^{k_n}$. If $n$ is non-ambiguous, let $\nabla G_\alpha(\theta) = (g_\alpha'(\theta(i)))_{i=1}^{k_n}$ and let $\nabla^2 G(\theta) = \mathrm{diag}\bigl(g_\alpha''(\theta(i))_{i=1}^{k_n}\bigr)$.

Then for some (random) vectors $\tilde{\tau}$ and $\tilde{\theta}$ with $\tilde{\tau}(i)$ (resp. $\tilde{\theta}(i)$) between $\tau(i)$ and $\theta_0(i)$ (resp. between $\widehat{\theta}_n(i)$ and $\theta_0(i)$) for all $i = 1, \ldots, k_n$:

$$\begin{aligned}\sqrt{n}\bigl(G_{n,\alpha}(T) - G_\alpha(\hat{\theta})\bigr) &= (H_n - \Delta_n(\theta_0))^T \nabla G_\alpha(\theta_0) \\ &\quad + \frac{1}{2\sqrt{n}} H_n^T \nabla^2 G_\alpha(\tilde{\tau}) H_n \\ &\quad - \frac{1}{2\sqrt{n}} \Delta_n^T(\theta_0) \nabla^2 G_\alpha(\tilde{\theta}) \Delta_n(\theta_0).\end{aligned}$$

This follows from

$$\begin{aligned}&\sqrt{n}\left(G_{n,\alpha}(T) - G_\alpha(\theta_0)\right) \\ &= -\sqrt{n} \sum_{i=k_n+1}^{+\infty} g_\alpha(\theta_0(i)) + \sqrt{n} \sum_{i=1}^{k_n} (g_\alpha(\tau(i)) - g_\alpha(\theta_0(i))) \\ &= -\sqrt{n} \sum_{i=k_n+1}^{+\infty} g_\alpha(\theta_0(i)) + H_n^T \nabla G_\alpha(\theta_0) + R_n\end{aligned}$$

with

$$R_n = \frac{1}{2\sqrt{n}} H_n^T \nabla^2 G_\alpha(\tilde{\tau}) H_n.$$

Meanwhile

$$\begin{aligned}&\sqrt{n}\left(G_{n,\alpha}(\widehat{\theta}_n) - G_\alpha(\theta_0)\right) \\ &= -\sqrt{n} \sum_{i=k_n+1}^{+\infty} g_\alpha(\theta_0(i)) + \sqrt{n} \sum_{i=1}^{k_n} \left(g_\alpha(\widehat{\theta}_n(i)) - g_\alpha(\theta_0(i))\right) \\ &= -\sqrt{n} \sum_{i=k_n+1}^{+\infty} g_\alpha(\theta_0(i)) + \Delta_n^T(\theta_0) \nabla G_\alpha(\theta_0) + \tilde{R}_n\end{aligned}$$

with

$$\tilde{R}_n = \frac{1}{2\sqrt{n}} \Delta_n^T(\theta_0) \nabla^2 G_\alpha(\tilde{\theta}) \Delta_n(\theta_0)$$

Now recall that $\gamma_{n,\alpha} = \operatorname{Var}\left((H_n - \Delta_n(\theta_0))^T \nabla g_\alpha(\theta_0)\right)$. Let $(\epsilon_n)_n$ be a sequence tending to 0 as $n$ tends to infinity.

For any interval $I$ of $\mathbb{R}$, let $I_{\epsilon_n}$ be the $\epsilon_n$-blowup of $I$: $I_{\epsilon_n} = \{x : \exists y \in I, |x - y| \leq \epsilon_n\}$. Then, for some positive constant $C$

$$\begin{aligned}&\sup_{I \in \mathcal{I}} \left| P_{H_n | X_{1:n}} \left( \frac{\sqrt{n}(G_{n,\alpha}(T) - G_{n,\alpha}(\widehat{\theta}_n))}{\gamma_{n,\alpha}} \in I \right) - \Phi(I) \right| \\ &\leq \sup_{I \in \mathcal{I}} \left| P_{H_n | X_{1:n}} \left( \frac{(H_n - \Delta_n(\theta_0))^T \nabla g_\alpha(\theta_0)}{\gamma_{n,\alpha}} \in I_{\epsilon_n} \right) - \Phi(I) \right| \\ &\quad + P_{H_n | X_{1:n}} \left( |R_n - \tilde{R}_n| \geq \epsilon_n C \right)\end{aligned}$$

The first summand on the right-hand-side is easily dealt with by applying the Bernstein-Von Mises Theorem:

$$\sup_{I \in \mathcal{I}} \left| P_{H_n | X_{1:n}} \left( \frac{(H_n - \Delta_n(\theta_0))^T \nabla g_\alpha(\theta_0)}{\gamma_{n,\alpha}} \in I_{\epsilon_n} \right) - \Phi(I) \right|$$
$$\leq \sup_{I \in \mathcal{I}} |\Phi(I_{\epsilon_n}) - \Phi(I)| + \left\| P_{H_n | X_{1:n}} - \mathcal{N}(\Delta_n(\theta_0), I^{-1}(\theta_0)) \right\|$$
$$\leq \frac{\epsilon_n}{\sqrt{2\pi}} + \left\| P_{H_n | X_{1:n}} - \mathcal{N}(\Delta_n(\theta_0), I^{-1}(\theta_0)) \right\| .$$

Now,

$$P_{H_n | X_{1:n}} \left( |R_n - \tilde{R}_n| \geq C \epsilon_n \right) \leq P_{H_n | X_{1:n}} \left( R_n \geq \frac{C \epsilon_n}{2} \right) + \mathbb{1}_{\tilde{R}_n \geq \frac{C \epsilon_n}{2}},$$

as $\tilde{R}_n$ is $X_{1:n}$-measurable. Theorem 3.12 follows if it is possible to choose a sequence $(\epsilon_n)_n$ such that both terms in the upper bound tend to 0 in $\mathbb{P}_{n,0}$ probability.

Let us focus for the moment on the first term. We aim at proving that the following upper-bound holds for large enough $n$ and some positive constant $D$:

$$P_{H_n | X_{1:n}} \left( |R_n| \geq \frac{C \epsilon_n}{2} \right)$$
$$\leq 5 P_{H_n | X_{1:n}} \left( H_n^T I(\theta_0) H_n \geq D \epsilon_n \sqrt{n \inf_{j \leq k_n} \theta_0(j)} \right). \quad (5.1)$$

As $g''_\alpha$ is monotone, the following inequalities hold:

$$|g''_\alpha(\tilde{\tau}(i))| \leq \max \{ |g''_\alpha(\tau(i))|; |g''_\alpha(\theta_0(i))| \} \leq |g''_\alpha(\tau(i))| + |g''_\alpha(\theta_0(i))|$$

This entails

$$\sqrt{n} R_n \leq H_n^T \nabla^2 G_\alpha(\theta_0) H_n + H_n^T \nabla^2 G_\alpha(\tilde{\tau}) H_n .$$

$$P_{H_n | X_{1:n}} (R_n \geq C \epsilon_n) \leq P_{H_n | X_{1:n}} \left( H_n^T \nabla^2 G_\alpha(\theta_0) H_n \geq \frac{C \epsilon_n \sqrt{n}}{2} \right)$$
$$+ P_{H_n | X_{1:n}} \left( H_n^T \nabla^2 G_\alpha(\tilde{\tau}) H_n \geq \frac{C \epsilon_n \sqrt{n}}{2} \right) .$$

Henceforth, let $C_\alpha = \frac{1}{\alpha(\alpha-1)}$ for $\alpha \neq 1$ and $C_1 = 1$. Note that for all positive $x$, $C_\alpha g''_\alpha(x) = x^{\alpha-2}$ and $C_\alpha H_n^T \nabla^2 G_\alpha(\theta_0) H_n = \sum_{i=1}^{k_n} \frac{H_n^2(i)}{\theta_0(i)} \theta_0(i)^{\alpha-1}$.

For $\alpha \geq 1$,

$$P_{H_n | X_{1:n}} \left( H_n^T \nabla^2 G_\alpha(\theta_0) H_n \geq \frac{C \epsilon_n \sqrt{n}}{2} \right) \leq P_{H_n | X_{1:n}} \left( H_n^T I(\theta_0) H_n \geq \frac{C_\alpha C \epsilon_n \sqrt{n}}{2} \right).$$

Meanwhile, for $\frac{1}{2} \leq \alpha < 1$, the obvious fact $\sup_i (\theta_0(i))^{\alpha-1} \leq (\inf_{j \leq k_n} \theta_0(j))^{-1/2}$, implies

$$\sum_{i=1}^{k_n} H_n^2(i) \theta_0(i)^{\alpha-2} \leq \frac{H_n^T I(\theta_0) H_n}{\sqrt{n \inf_{j \leq k_n} \theta_0(j)}},$$

so, we get

$$P_{H_n|X_{1:n}}\left(\sum_{i=1}^{k_n} H_n^2(i)\theta_0(i)^{\alpha-2} \geq C_\alpha C\epsilon_n\sqrt{n}/2\right)$$
$$\leq P_{H_n|X_{1:n}}\left(H_n^T I(\theta_0)H_n \geq \frac{C_\alpha C\epsilon_n}{2}\sqrt{n\inf_{j\leq k_n}\theta_0(j)}\right).$$

On the other hand,

$$\sum_{i=1}^{k_n} H_n^2(i)\tau(i)^{\alpha-2}$$
$$\leq \sum_{i=1}^{k_n} \frac{H_n^2(i)}{\theta_0(i)}\theta_0(i)^{\alpha-2}\left(1+\sqrt{n}\frac{\sup_{i\leq k_n}\frac{|H_n(i)|}{\sqrt{\theta_0(i)}}}{\sqrt{n\inf_{j\leq k_n}\theta_0(j)}}\right)^{\alpha-2}$$
$$\leq \sum_{i=1}^{k_n} \frac{H_n^2(i)}{\theta_0(i)}\theta_0(i)^{\alpha-2}\left(1+\frac{\left(\sum_{i\leq k_n}\frac{H_n^2(i)}{\theta_0(i)}\right)^{1/2}}{\sqrt{n\inf_{j\leq k_n}\theta_0(j)}}\right)^{\alpha-2}.$$

Hence,

$$P_{H_n|X_{1:n}}\left(\sum_{i=1}^{k_n}(H^n(i))^2\tau(i)^{\alpha-2} \geq \frac{C_\alpha C\epsilon_n}{2\sqrt{n}}\right)$$
$$\leq P_{H_n|X_{1:n}}\left(\sum_{i=1}^{k_n}\frac{H_n^2(i)}{\theta_0(i)}\theta_0(i)^{\alpha-2} \geq 2^{\alpha-3}C_\alpha C\epsilon_n\sqrt{n}\right)$$
$$+ P_{H_n|X_{1:n}}\left(\sum_{i\leq k_n}\frac{H_n^2(i)}{\theta_0(i)} \geq n\inf_{j\leq k_n}\theta_0(j)\right).$$

We may now sum up those inequalities:

$$P_{H_n|X_{1:n}}\left(|R_n|\geq \frac{C\epsilon_n}{2}\right) \leq \sum_{i=1}^{5} P_{H_n|X_{1:n}}\left(H_n^T I(\theta_0)H_n \geq A_{i,n}\right)$$

with

$$\begin{aligned}
A_{1,n} &= \frac{C_\alpha C}{2}\epsilon_n\sqrt{n}\\
A_{2,n} &= \frac{C_\alpha C}{2}\sqrt{n\inf_{j\leq k_n}\theta_0(j)}\epsilon_n\\
A_{3,n} &= 2^{\alpha-2}A_{1,n}\\
A_{4,n} &= 2^{\alpha-2}A_{2,n}\\
A_{5,n} &= n\inf_{j\leq k_n}\theta_0(j).
\end{aligned}$$

This is enough to prove Inequality (5.1). Up to $\|P_{H_n|X_{1:n}} - \mathcal{N}_{k_n}\|$, the posterior probability of the event $H_n^T I(\theta_0) H_n \geq u_n$ equals $\mathcal{N}_{k_n}\left\{H_n^T I(\theta_0) H_n \geq u_n\right\}$, that is the probability that a non-central $\chi_{k_n}^2$-distributed random variable with non-centrality parameter $\Delta_n^T(\theta_0) I_n(\theta_0) \Delta_n(\theta_0) = V_n(\theta_0)^2$ exceeds $u_n$. The latter probability is upper-bounded by

$$\mathbb{1}_{V_n^2(\theta_0) \geq \frac{u_n}{4}} + P\left(\xi_n^2 \geq \frac{u_n}{4}\right)$$

where $\xi_n^2$ follows a $\chi_{k_n}^2$ distribution. The probability that $\xi_n^2$ is larger than $u_n/4$ may be upper bounded using Cirelson's inequality. As soon as $(\epsilon_n)_n$ satisfies

$$\epsilon_n \frac{\sqrt{n \inf_{j \leq k_n} \theta_0(j)}}{k_n} \to +\infty,$$

this probability tends to 0.

The $\mathbb{P}_{n,0}$-probability that the non-centrality parameter $V_n^2(\theta_0)$ is large may be upper bounded using Inequality (4.3), and invoking the Bernstein-Von Mises Theorem to handle $\|P_{H_n|X_{1:n}} - \mathcal{N}(\Delta_n(\theta_0), I^{-1}(\theta_0))\|$,

$$P_{H_n|X_{1:n}}\left(|R_n| \geq \frac{C\epsilon_n}{2}\right) \to 0$$

in $\mathbb{P}_{n,0}$-probability.

Using the same approach as before, one establishes that for large enough $n$ and some positive constant $D$:

$$\mathbb{P}_{n,0}\left(|\tilde{R}_n| \geq \frac{C\epsilon_n}{2}\right) \leq 5\,\mathbb{P}_{n,0}\left(V_n^2(\theta_0) \geq D\epsilon_n \sqrt{n \inf_{j \leq k_n} \theta_0(j)}\right).$$

The right-hand-side may be upper-bounded using again (4.3). It tends to 0 as soon as $\epsilon_n \sqrt{n \inf_{j \leq k_n} \theta_0(j)}/k_n \to +\infty$.

## 6. Proof of Theorem 3.13

We have already proved that $\tilde{R}_n = o_{\mathbb{P}_{n,0}}(1)$, so that under the assumptions of Theorem 3.13, equation (5.1) translates into

$$\sqrt{n}\left(G_{n,\alpha}(\widehat{\theta}_n) - G_\alpha(\theta_0)\right) = o(1) + \Delta_n^T(\theta_0)\nabla G_\alpha(\theta_0) + o_{\mathbb{P}_{n,0}}(1)$$

and the result follows from Berry-Essen Theorem.

## Appendix A: Contiguity

We first prove Lemma 2.3.

*Proof.* Let us first notice that if $\sigma_n^2(h_n)$ tends toward $\sigma^2 > 0$, there exists some $M > 0$ such that $h_n \in \mathcal{E}(M)$. The contiguity proof follows from a straightforward analysis of the log-likelihood ratio and an invocation of Le Cam's first Lemma (van der Vaart, 2002).

A Taylor expansion of the logarithm leads to

$$\log \frac{\mathbb{P}_{n,h_n}}{\mathbb{P}_{n,0}}(\mathbf{x})$$
$$= \sum_{i=0}^{k_n} N_i \log\left(1 + \frac{h(i)}{\sqrt{n}\theta_0(i)}\right)$$
$$= \frac{1}{\sqrt{n}} \sum_{i=0}^{k_n} N_i \frac{h(i)}{\theta_0(i)} - \frac{1}{2n} \sum_{i=0}^{k_n} N_i \frac{h_n(i)^2}{\theta_0(i)^2} + \frac{1}{n} \sum_{i=0}^{k_n} N_i \frac{h_n(i)^2}{\theta_0(i)^2} R\left(\frac{h_n(i)}{\sqrt{n}\theta_0(i)}\right).$$

The proof consists in checking the three following points:

1. the remainder term converges in probability toward 0.
2. the first summand converges in distribution toward $\mathcal{N}(0, \sigma^2)$.
3. the middle term converges in probability toward $-\sigma^2/2$.

Let us check the first point. As a matter of fact:

$$\mathbb{P}_{n,\theta_0}\left(\frac{1}{n}\sum_{i=0}^{k_n} N_i \frac{h_n(i)^2}{\theta_0(i)^2} R\left(\frac{h_n(i)}{\sqrt{n}\theta_0(i)}\right)\right) \leq M \cdot \sup_{i \leq k_n} \left| R\left(\frac{h_n(i)}{\sqrt{n}\theta_0(i)}\right)\right|$$

But
$$\sup_{i \leq k_n} \frac{|h_n(i)|}{\sqrt{n}\theta_0(i)} \leq \frac{\sigma_n(h_n)}{\sqrt{n \inf_{i \leq k_n} \theta_0(i)}} = o(1).$$

In order to check the the second point, note that the random variable $Z_n(h_n) = \frac{1}{\sqrt{n}} \sum_{i=0}^{k_n} N_i \frac{h(i)}{\theta_0(i)}$ can be rewritten as a sum of i.i.d. random variables:

$$Z_n(h_n) = \frac{1}{\sqrt{n}} \sum_{j=1}^{n} Y_j$$

with
$$Y_j = \sum_{i=0}^{k_n} \left(\frac{\mathbf{1}_{X_j = i} - \theta_0(i)}{\theta_0(i)}\right) h_n(i);$$

Under $\mathbb{P}_{n,0}$, each random variable $Y_j$ is equal to $\frac{h_n(i)}{\theta_0(i)}$ with probability $\theta_0(i)$, so that
$$\mathbb{P}|Y_j|^3 = \sum_{i=0}^{k_n} \frac{|h_n(i)|^3}{\theta_0(i)^2}.$$

$$\sum_{i=0}^{k_n} \frac{|h_n(i)|^3}{\theta_0(i)^2} \leq \left(\sup_{i \leq k_n} \frac{h_n(i)}{\sqrt{n}\theta_0(i)}\right) \sqrt{n}\sigma_n^2(h_n),$$

that is
$$\sum_{i=0}^{k_n} \frac{|h(i)|^3}{\theta_0(i)^2} = o\left[\sqrt{n}\sigma_n^3(h_n)\right]$$

as $\sigma_n^2(h_n)$ is bounded and bounded away from 0. The Berry-Essen Theorem (Dudley, 2002) entails that as $n$ tends to infinity, $Z_n(h_n)$ converges in distribution toward $\mathcal{N}(0,\sigma^2)$.

Finally, the middle term $\frac{1}{2n}\sum_{i=0}^{k_n} N_i \frac{h_n(i)^2}{\theta_0(i)^2}$ converges in probability toward $\frac{1}{2}\sigma_n(h_n)^2$. Indeed, let

$$U_n(h_n) = \frac{1}{n}\sum_{i=0}^{k_n} N_i \frac{h_n(i)^2}{\theta_0(i)^2} - \sigma_n^2(h_n) = \frac{1}{n}\sum_{j=1}^{n}(\xi_j(h_n) - E_0(\xi_j(h_n)))$$

with
$$\xi_j(h_n) = \sum_{i=0}^{k_n} 1_{X_j=i}\frac{h_n(i)^2}{\theta_0(i)^2}.$$

Then
$$\begin{aligned}\operatorname{Var}(U_n(h_n)) &= \frac{1}{n}\operatorname{Var}(\xi_1(h_n)) \\ &\leq \frac{1}{n}\sum_{i=0}^{k_n}\frac{h_n(i)^4}{\theta_0(i)^3} \\ &\leq \frac{M^2}{n\inf_{i\leq k_n}\theta_0(i)}.\end{aligned}$$

Hence, the sequence of distributions of likelihood ratios $\frac{\mathbb{P}_{n,h}}{\mathbb{P}_{n,0}}(X_{1:n})$ converges weakly toward a log-normal distribution with parameters $-\sigma^2/2$ and $\sigma^2$. The Lemma follows directly from Le Cam's first Lemma (van der Vaart, 1998). □

**Lemma A.1.** *Let $\theta_0$ denote a probability mass function over $\mathbb{N}_*$. If the sequence of truncation levels $(k_n)_{n\in\mathbb{N}}$ satisfy Condition 2.1, the sequences $(\mathbb{P}_{n,0})_n$ and $(\mathbb{P}_{n,A})_n$ are mutually contiguous.*

*Proof.* Let $(B_n)$ be a sequence of events where for each $n$, $B_n \subseteq \{0,\ldots,k_n\}^n$. Then
$$\mathbb{P}_{n,A}(B_n) \leq \sup_{\sigma_n(h)^2 \leq A} \mathbb{P}_{n,h}(B_n)$$
so that for some sequence $(h_n)_n$ such that for all $n$, $\sigma_n^2(h_n) \leq A$,
$$\limsup_{n\to+\infty}\mathbb{P}_{n,A}(B_n) \leq \limsup_{n\to+\infty}\mathbb{P}_{n,h_n}(B_n).$$

But as $(\theta_0(i))$ decreases to 0 at infinity, $\mathcal{E}_{\theta_0,k_n}(A)$ is a finite dimensional closed subset of a compact set in $\ell^2(\mathbb{N}_*)$, so that one may extract a subsequence $(h_{n_p})_p$ such that $\sigma^2_{n_p}(h_{n_p}) \to \sigma^2$ for some $\sigma^2$ as $p \to +\infty$, and such that

$$\limsup_{n\to+\infty} \mathbb{P}_{n,A}(B_n) \leq \lim_{p\to+\infty} \mathbb{P}_{n_p,h_{n_p}}(B_{n_p}).$$

Applying Lemma 2.3 gives that $\mathbb{P}_{n,0}(B_n) \to 0$ implies $\mathbb{P}_{n,A}(B_n) \to 0$.
The reverse implication may be proved with the same reasoning using that

$$\inf_{\sigma^2_n(h) \leq A} \mathbb{P}_{n,h}(B_n) \leq \mathbb{P}_{n,A}(B_n).$$

□

## Appendix B: Distance in variation and conditioning

The obvious proof of the following folklore lemma is omitted.

**Lemma B.1.** *Let $P$ denote a probability distribution on some space $(\Omega, \mathcal{F})$. Let $A$ denote an event with non-null $P$-probability and let $P^A$ the conditional probability given $A$, that is $P^A(B) = P(A \cap B)/P(A)$ then*

$$\|P^A - P\| = P(A^c).$$

## Appendix C: Proof of inequality (3.9)

*Proof.*

$$\left\| \mathcal{N}^{M_n}_{k_n} - P^{M_n}_{H_n|X_{1:n}} \right\|$$

$$= \int \left(1 - \frac{d\mathcal{N}^{M_n}_{k_n}(h)}{dP^{M_n}_{H_n|X_{1:n}}(h)}\right)_+ dP^{M_n}_{H_n|X_{1:n}}(h)$$

$$= \int \left(1 - \frac{d\mathcal{N}^{M_n}_{k_n}(h) \int \frac{w_n(\theta_0 + \frac{g}{\sqrt{n}})}{d\mathcal{N}^{M_n}_{k_n}(g)} \frac{dP_{n,g}(X_{1:n})}{dP_{n,0}(X_{1:n})} d\mathcal{N}^{M_n}_{k_n}(g)}{w_n(\theta_0 + \frac{h}{\sqrt{n}}) \frac{dP_{n,h}(X_{1:n})}{dP_{n,0}(X_{1:n})}}\right)_+ dP^{M_n}_{H_n|X_{1:n}}(h).$$

Using the convexity of $x \mapsto (1-x)_+$ and Jensen inequality, the right-hand-side can be upper-bounded:

$$\left\| \mathcal{N}^{M_n}_{k_n} - P^{M_n}_{H_n|X_{1:n}} \right\|$$

$$\leq \iint \left(1 - \frac{d\mathcal{N}^{M_n}_{k_n}(h) w_n(\theta_0 + \frac{g}{\sqrt{n}}) \frac{dP_{n,g}(X_{1:n})}{dP_{n,0}(X_{1:n})}}{d\mathcal{N}^{M_n}_{k_n}(g) w_n(\theta_0 + \frac{h}{\sqrt{n}}) \frac{dP_{n,h}(X_{1:n})}{dP_{n,0}(X_{1:n})}}\right)_+ d\mathcal{N}^{M_n}_{k_n}(g) dP^{M_n}_{H_n|X_{1:n}}(h).$$

Now, the quadratic Taylor expansion of the log-likelihood ratio translates into

$$\frac{\mathbb{P}_{n,g}(X_{1:n}) d\mathcal{N}^{M_n}_{k_n}(h)}{\mathbb{P}_{n,h}(X_{1:n}) d\mathcal{N}^{M_n}_{k_n}(g)} = e^{\left\{\frac{1}{2}(A_n(g) - A_n(h)) + (C_n(g) - C_n(h))\right\}}.$$

Plugging this expansion into the upper-bound on $\text{NV}(M_n)$ leads to (3.9).

□

**Appendix D: Tail bounds for non-centered Pearson statistics**

This section provides a proof of Inequality (4.4). Recall from Section 4.2, that

$$\mathbb{P}_{n,h}\left(V_n(\theta_0) \leq s_n\right) \leq \mathbb{P}_{n,h}\left(\sum_{i=0}^{k_n} \mathbf{a}_i^* \frac{N_i - n\theta(i)}{\sqrt{n\theta_0(i)}} \leq -\sigma_n(h) + s_n\right).$$

where $\mathbf{a}_i^* = \frac{h(i)}{\sqrt{\theta_0(i)\sigma^2(h)}}$ for all $i \leq k_n$. Despite $\sum_{i=0}^{k_n} \mathbf{a}_i^* \frac{N_i - n\theta(i)}{\sqrt{n\theta_0(i)}}$ is just a sum of i.i.d. random variables, we found no obvious way to use classical exponential inequalities (either Hoeffding or Bernstein inequalities) to prove the tail bounds we need. Before resorting to classical inequalities, we split the sum into two pieces according to the signs of the coefficients $\mathbf{a}_i^*$. The two pieces are handled using negative association arguments.

Let $\mathcal{J} = \{i : i \leq k_n, \mathbf{a}_i^* \geq 0\}$ and $\mathcal{J}^c = \{i : i \leq k_n, \mathbf{a}_i^* < 0\}$. Note first that

$$\begin{aligned}
&\mathbb{P}_{n,h}\left(\sum_{i=0}^{k_n} \mathbf{a}_i^* \frac{N_i - n\theta(i)}{\sqrt{n\theta_0(i)}} \leq -\sigma_n(h) + s_n\right) \\
&\leq \mathbb{P}_{n,h}\left(\sum_{i \in \mathcal{J}} \mathbf{a}_i^* \frac{N_i - n\theta(i)}{\sqrt{n\theta_0(i)}} \leq -\frac{1}{2}(\sigma_n(h) - s_n)\right) \\
&\quad + \mathbb{P}_{n,h}\left(\sum_{i \in \mathcal{J}^c} \mathbf{a}_i^* \frac{N_i - n\theta(i)}{\sqrt{n\theta_0(i)}} \leq -\frac{1}{2}(\sigma_n(h) - s_n)\right) \\
&\leq \inf_{\lambda<0} \exp\left(\log \mathbb{P}_{n,h}\left[\exp\left(\sum_{i \in \mathcal{J}} \lambda \mathbf{a}_i^* \frac{N_i - n\theta(i)}{\sqrt{n\theta_0(i)}}\right)\right] + \frac{\lambda}{2}(\sigma_n(h) - s_n)\right) \\
&\quad + \inf_{\lambda>0} \exp\left(\log \mathbb{P}_{n,h}\left[\exp\left(-\sum_{i \in \mathcal{J}^c} \lambda \mathbf{a}_i^* \frac{N_i - n\theta(i)}{\sqrt{n\theta_0(i)}}\right)\right] - \frac{\lambda}{2}(\sigma_n(h) - s_n)\right).
\end{aligned}$$

Following Dubhashi and Ranjan (1998), a collection of random variables $Z_1, \ldots, Z_n$ is said to be negatively associated if for any $\mathcal{I} \subseteq \{1, \ldots, n\}$, for any functions $f : \mathbb{R}^{|\mathcal{I}|} \to \mathbb{R}$ and $g : \mathbb{R}^{\mathcal{I}^c} \to \mathbb{R}$ that are either both non-decreasing or both non-increasing,

$$\mathbb{P}\left[f(X_i : i \in \mathcal{I})g(X_i : i \in \mathcal{I}^c)\right] \leq \mathbb{P}\left[f(X_i : i \in \mathcal{I})\right]\mathbb{P}\left[g(X_i : i \in \mathcal{I}^c)\right].$$

By Theorem 14 from Dubhashi and Ranjan (1998), both sets of random variables $\left(\mathbf{a}_i^*(N_i - n\theta(i))/\sqrt{n\theta_0(i)}\right), i \in \mathcal{J}$ and $\left(\mathbf{a}_i^*(N_i - n\theta(i))/\sqrt{n\theta_0(i)}\right), i \in \mathcal{J}^c$ are negatively associated in the sense of Dubhashi and Ranjan (1998).

The logarithmic moment generating function of $\sum_{i \in \mathcal{I}} \mathbf{a}_i^* \frac{N_i - n\theta(i)}{\sqrt{n\theta_0(i)}}$ satisfies

$$\log \mathbb{P}_{n,h} e^{\lambda \sum_{i \in \mathcal{I}} \mathbf{a}_i^* \frac{N_i - n\theta(i)}{\sqrt{n\theta_0(i)}}} \leq \sum_{i \in \mathcal{I}} \log \mathbb{P}_{n,h} e^{\lambda \mathbf{a}_i^* \frac{N_i - n\theta(i)}{\sqrt{n\theta_0(i)}}}$$

where $\mathcal{I} = \mathcal{J}$ or $\mathcal{J}^c$. Each $N_i$ is binomially distributed with parameter $n$ and $\theta_i$.

For $i \in \mathcal{J}, \mathbf{a}_i^* \geq 0$ so that for $\lambda \leq 0$:

$$\log \mathbb{P}_{n,h} e^{\lambda \sum_{i \in \mathcal{J}} \mathbf{a}_i^* \frac{(N_i - n\theta(i))}{\sqrt{n\theta_0(i)}}} \leq \sum_{i \in \mathcal{J}} \frac{(\lambda \mathbf{a}_i^*)^2 \theta(i)}{2\theta_0(i)}.$$

Note that

$$\begin{aligned}
\sum_{i \in \mathcal{J}} \mathbf{a}_i^{*2} \frac{\theta(i)}{\theta_0(i)} &\leq \sum_{i \in \mathcal{J}} \mathbf{a}_i^{*2} + \sum_{i \in \mathcal{J}} \mathbf{a}_i^{*2} \frac{1}{\sqrt{n\theta_0(i)}} \frac{h(i)}{\sqrt{\theta_0(i)}} \\
&\leq \sum_{i \in \mathcal{J}} \mathbf{a}_i^{*2} + \frac{\left(\sum_{i \in \mathcal{J}} \mathbf{a}_i^{*4}\right)^{1/2}}{\sqrt{n \inf_{i \leq k_n} \theta_0(i)}} \left(\sum_{i \in \mathcal{J}} \frac{h^2(i)}{\theta_0(i)}\right)^{1/2} \\
&\leq 1 + \frac{\sigma_n(h)}{\sqrt{n \inf_{i \leq k_n} \theta_0(i)}}.
\end{aligned}$$

Hence

$$\inf_{\lambda < 0} \exp\left(\log \mathbb{P}_{n,h}\left[\exp\left(\sum_{i \in \mathcal{J}} \lambda \mathbf{a}_i^* \frac{N_i - n\theta(i)}{\sqrt{n\theta_0(i)}}\right)\right] + \frac{\lambda}{2}(\sigma_n(h) - s_n)\right)$$

$$\leq \exp\left(-\frac{(\sigma_n(h) - s_n)^2}{8\left(1 + \frac{\sigma_n(h)}{\sqrt{n \inf_{i \leq k_n} \theta_0(i)}}\right)}\right).$$

If $\sigma_n(h) \geq 2s_n$, $\sigma_n(h) - s_n \geq \sigma_n(h)/2$ and the last term may be upper-bounded by $\exp\left(-\frac{1}{64}\left[\sigma_n^2(h) \wedge \sigma_n(h)\sqrt{n \inf_{i \leq k_n} \theta_0(i)}\right]\right)$.

Now, if $i \in \mathcal{J}^c$, $h(i) \leq 0$ so that $h(i) \geq -\sqrt{n}\theta_0(i)$ which entails $-\mathbf{a}_i^*/\sqrt{n\theta_0(i)} \leq \frac{1}{\sigma_n(h)}$. For any $\lambda \geq 0$

$$\sum_{i \in \mathcal{J}^c} \log \mathbb{P}_{n,h}\left[\exp\left(-\lambda \mathbf{a}_i^* \frac{N_i - n\theta(i)}{\sqrt{n\theta_0(i)}}\right)\right] \leq \frac{\sum_{i \in \mathcal{J}^c}(\mathbf{a}_i^*)^2 \theta(i) \lambda^2}{2\theta_0(i)\left(1 - \frac{\lambda}{\sigma_n(h)}\right)}$$

$$\leq \lambda^2 \frac{\left(1 + \frac{\sigma_n(h)}{\sqrt{n \inf_{i \leq k_n} \theta_0(i)}}\right)}{2\left(1 - \frac{\lambda}{\sigma_n(h)}\right)}.$$

Hence

$$\inf_{\lambda > 0} \exp\left(\log \mathbb{P}_{n,h}\left[\exp\left(-\sum_{i \in \mathcal{J}^c} \lambda \mathbf{a}_i^* \frac{N_i - n\theta(i)}{\sqrt{n\theta_0(i)}}\right)\right] - \frac{\lambda}{2}(\sigma_n(h) - s_n)\right)$$

$$\leq \exp\left(-\frac{(\sigma_n(h) - s_n)^2}{8\left(1 + \frac{\sigma_n(h)}{\sqrt{n \inf_{i \leq k_n} \theta_0(i)}} + \frac{\sigma_n(h) - s_n}{\sigma_n(h)}\right)}\right).$$

If $\sigma_n(h) \geq 2s_n$, $\sigma_n(h) - s_n \geq \sigma_n(h)/2$ and the right-hand-side is upper-bounded by

$$\exp\left(-\frac{\left(\sigma_n^2(h) \wedge \sqrt{n \inf_{i \leq k_n} \theta_0(i)}\sigma_n(h)\right)}{96}\right).$$